\documentclass[onefignum,onetabnum]{siamart171218}

\usepackage{graphicx}
\usepackage{tikz}
\usepackage{booktabs}
\usepackage{amsmath}
\usepackage{bm}
\theoremstyle{plain}
\newtheorem{example}{Example}[section]

\usepackage[top=1in, bottom=1in, left=1in, right=1in]{geometry}

\newcommand{\bd}[1]{\boldsymbol{#1}}

\newcommand{\wh}[1]{\widehat{#1}}
\newcommand{\wt}[1]{\widetilde{#1}}

\newcommand{\veps}{\varepsilon}
\newcommand{\eps}{\epsilon}

\newcommand{\rd}{\,\mathrm{d}}

\newcommand{\comment}[1]{}

\newcommand{\hy}[1]{\textcolor{blue}{HY: #1}}

\usepackage{lipsum}
\usepackage{amsfonts}
\usepackage{graphicx}
\usepackage{epstopdf}
\usepackage{algorithmic}
\ifpdf
  \DeclareGraphicsExtensions{.eps,.pdf,.png,.jpg}
\else
  \DeclareGraphicsExtensions{.eps}
\fi

\newsiamremark{remark}{Remark}
\newsiamremark{hypothesis}{Hypothesis}
\crefname{hypothesis}{Hypothesis}{Hypotheses}
\newsiamthm{claim}{Claim}

\headers{Multiscale and Nonlocal Learning for PDEs}{D. Ricardo, H. Jingwei, and Y. Haizhao}

\title{Multiscale and Nonlocal Learning for PDEs using Densely Connected RNNs\thanks{Submitted to the editors DATE.
\funding{RD was partially supported by the Fog Research Institute under contract no.~FRI-454.  JH was partially supported by NSF CAREER grant DMS-1654152. HY was partially supported by NSF CAREER grant DMS-1945029.
}}}

\author{Ricardo A. Delgadillo\thanks{Department of Civil and Environmental Engineering, National University of Singapore, Singapore (\email{ceerad@nus.edu.sg}).}
\and Jingwei Hu\thanks{Department of Mathematics, Purdue University, West Lafayette, IN 47907, USA (\email{jingweihu@purdue.edu}).}
\and Haizhao Yang\thanks{Department of Mathematics, Purdue University, West Lafayette, IN 47907, USA (\email{haizhao@purdue.edu}).}}

\usepackage{amsopn}

\ifpdf
\hypersetup{
  pdftitle={Multiscale and Nonlocal Learning for PDEs using Densely Connected RNNs},
  pdfauthor={D. Ricardo, H. Jingwei, and Y. Haizhao}
}
\fi

\begin{document}
\maketitle
\begin{abstract}
Learning time-dependent partial differential equations (PDEs) that govern evolutionary observations is one of the core challenges for data-driven inference in many fields. In this work, we propose to capture the essential dynamics of numerically challenging PDEs arising in multiscale modeling and simulation -- kinetic equations. These equations are usually nonlocal and contain scales/parameters that vary by several orders of magnitude. We introduce an efficient framework, Densely Connected Recurrent Neural Networks (DC-RNNs), by incorporating a multiscale ansatz and high-order implicit-explicit (IMEX) schemes into RNN structure design to identify analytic representations of multiscale and nonlocal PDEs from discrete-time observations generated from heterogeneous experiments.  If present in the observed data, our DC-RNN can capture transport operators, nonlocal projection or collision operators, macroscopic diffusion limit, and other dynamics. We provide numerical results to  demonstrate the advantage of our proposed framework and compare with existing methods.
\end{abstract}

\begin{keywords}
 Multiscale; Nonlocal; Kinetic Equation; Time-Dependent PDE Recovery; Machine Learning; Densely Connected Recurrent Neural Network.
\end{keywords}

\begin{AMS}
  35C20; 35C99; 68T99.
\end{AMS}

\section{Introduction}

Data-driven discovery of PDEs is experiencing unprecedented development over the past few years, wherein various kinds of PDEs (featuring, for example, time dependence and nonlinearity) have been studied. In this work, we consider the learning problem for a class of PDEs that involve multiple time/spatial scales and nonlocal operators -- kinetic equations. These are an important class of equations in multiscale modeling hieriarchy which bridges microscopic atomistic models (such as N-body Newton equations) and macroscopic continuum models (such as Navier-Stokes equations). For a variety of scientific problems ranging from gas/plasma dynamics, radiative transfer to social/biological systems, kinetic equations have demonstrated their ability to accurately model the dynamics of complex systems \cite{Villani02}. To the best of our knowledge, learning of multiscale kinetic equations, albeit important, has never been explored in the literature.


Specifically, we are interested in developing an efficient symbolic neural network to fit time-dependent data for a class of multiscale kinetic equations. The overall goal is to identify an explicit formula of the map $\mathcal{F}$ that determines the evolution $u(\bd{x},t)\rightarrow u(\bd{x},t+\Delta t)$ for $\bd{x}\in\Omega$ and $\Delta t>0$. Therefore, a symbolic neural network $\mathcal{F}(u;\bd{\theta},w_{\veps})$ with parameters $\bd{\theta}$ and $w_{\veps}$ is constructed and the following loss function is minimized to find the best parameter set:
\begin{equation}\label{eq:intro}
L(\bd{\theta},w_{\veps})=\dfrac{1}{N_{t}}\sum_{j=1}^{N_{t}}\left\|u(\bd{x},t_{j+1})-u(\bd{x},t_{j})-\int_{t_{j}}^{t_{j}+\Delta t}\mathcal{F}(u(\bd{x},s);\bd{\theta},w_{\veps})\,ds\right\|_{L^1(\Omega)}.
\end{equation}
Note that $\mathcal{F}$ approaches the correct model as $L(\bd{\theta},w_{\veps})\rightarrow 0$. The choice of the norm above is flexible. In this paper, we focus on the $L^1$-norm because our numerical experiments show that it is slightly better than others, e.g., the $L^2$-norm. Due to the multiscale and nonlocal feature of our target equations, existing learning schemes may not be efficient. We will propose novel symbolic neural networks, new formulations of the loss function in \eqref{eq:intro}, and new regularization methods in this paper to tackle this challenge.

Our first main contribution is a new symbolic neural network $\mathcal{F}(u;\bd{\theta},w_{\veps})$ built with multiscale and nonlocal features. The key idea for capturing multiscale phenomena is to construct $\mathcal{F}$ as a sum of different components at different scales  of order $\veps_{pred}^{n}$, where $n$ is an integer degree and $\veps_{pred}$ is a trainable multiscale separator defined by:
\begin{equation}\label{eqn:eps}
\veps_{pred}(w_{\veps})=\dfrac{1}{2}(\text{tanh}(w_{\veps})+1)
\end{equation}
with $w_{\veps}$ as a trainable parameter. In particular, we propose

\begin{equation}\label{eq:intro2}
\mathcal{F}(u;\bd{\theta},w_{\veps})=\sum_{n=0}^{N}\dfrac{1}{\veps_{pred}^{n}(w_{\veps})}\mathcal{F}^{n}(u;\bd{\theta}_{n}),
\end{equation}
where $\bd{\theta}:=(\bd{\theta}_{1},\bd{\theta}_{2},\cdots,\bd{\theta}_{n})$ and $\mathcal{F}^{n}(u;\bd{\theta}_{n})$ is the network at the $n$-th scale. Thus, unlike conventional deep learning recovery algorithms as in \cite{DBLP:journals/corr/abs-1711-10561,DBLP:journals/corr/abs-1711-10566,BinD18,harlim2020machine,WU2020109307,zslg18}, our algorithm is aware of different scales and thus more accurately captures different components at scale $\mathcal{O}(\veps_{pred}^{n})$. 

The key idea to make $\mathcal{F}(u;\bd{\theta},w_{\veps})$ capable of capturing nonlocal phenomena is to incorporate nonlocal operators in $\mathcal{F}^{n}$ in \eqref{eq:intro2} to construct $\mathcal{F}$. Conventionally, $\mathcal{F}$ is typically constructed as a linear combination of mathematical operators in a pre-specified dictionary, and the combination coefficients are learned via minimizing \eqref{eq:intro} with sparsity regularization to obtain sparse linear combinations as in \cite{Kaiser2018SparseIO,Schaeffer2017ExtractingSH,Mangan2017ModelSF,Brunton3932,zslg18}. For high-dimensional problems, constructing such a dictionary can be very costly. Hence, we will apply symbolic recurring neural network (RNN) of mathematical operators as in \cite{BinD17,BinD18} without specifying a large dictionary. Intuitively, due to the high expressiveness of our symbolic RNNs, the class of RNNs with different parameters can form a large dictionary without pre-specifying a costly dictionary. It might be computationally more efficient to use symbolic RNNs to classify the dynamics of data and choose a trainable symbolic model to model data.

The most basic elements in our RNN are a set of (either local or nonlocal) mathematical operators $\mathcal{A}_{1},\cdots,\mathcal{A}_{n}$ modeling the dynamics in kinetic equations, such as transport, collision, and diffusion operators. The trainable compositions of these basic operators form a basis of our RNN, i.e., each term $\mathcal{F}^{n}$ in \eqref{eq:intro2} is a trainable linear combination of the compositions defined below:
\begin{equation}\label{eq:rnnterms}
\mathcal{A}_{\pi(1)}\circ\cdots\circ\mathcal{A}_{\pi(m)},
\end{equation}
where $\pi=\begin{pmatrix}\pi(1),\cdots,\pi(m)\end{pmatrix}\in \mathbb{Z}^m$ with entries in $\{1,\cdots, n\}$. More precisely, we have
\begin{equation}\label{eq:rnnterms2}
\mathcal{F}^{n}(u;\bd{\theta}_{n})=\sum_{m\geq 1}\sum_{\pi\in\mathcal{D}}a_{\pi(1),\cdots,\pi(m)}(\bd{\theta}_{n})\mathcal{A}_{\pi(1)}\circ\cdots\circ\mathcal{A}_{\pi(m)}(u),
\end{equation}
where coefficients $a_{\pi(1),\cdots,\pi(m)}(\bd{\theta}_{n})$ depend on trainable parameters $\bd{\theta}_{n}$, and $\mathcal{D}$ is a set of index vectors $\pi$ specified by our symbolic RNN as we shall see later. Similar to polynomial regression \cite{doi:10.1080/00224065.1997.11979762,GERGONNE1974439}, our RNN returns a multivariate polynomial of the operators $\mathcal{A}_{1},\cdots,\mathcal{A}_{n}$. Due to the expressive power of neural networks \cite{yarotsky2017,Shen2,Shen3,li2020curse,Poggio2017WhyAW,liang2021reproducing}, our symbolic RNN of a small size can generate a sufficiently large index vector set $\mathcal{D}$. The formulation in \eqref{eq:rnnterms2} is also natural in physics, equations derived from asymptotic analysis often have recursive structure similar to the compositional operators in \eqref{eq:rnnterms2}, e.g., see \cite{SR09}. 

Our second main contribution is to propose novel loss functions based on high-order IMEX schemes to discretize of the integral in \eqref{eq:intro}. The most typical numerical method, the forward-Euler scheme, results in the loss function:
\begin{equation}\label{eq:intro3}
L(\bd{\theta},w_{\veps})=\dfrac{1}{N_{t}}\sum_{j=1}^{N_{t}}\left\|u(\bd{x},t_{j+1})-u(\bd{x},t_{j})-\Delta t \mathcal{F}(u(\bd{x},t_{j});\bd{\theta},w_{\veps})\right\|_{L^1(\Omega)},
\end{equation}
which is commonly used in the discovery of governing equations. Though higher order approximations using multistep methods have been investigated in \cite{DBLP:journals/corr/abs-1711-10561,keller2020discovery,Raissi2018MultistepNN,du2021discovery}, there is no existing research on the effectiveness of IMEX schemes in the literature of discovering governing equations. For kinetic equations, since they often contain non-stiff terms as well as stiff terms, the IMEX schemes are the natural choices and have demonstrated their power in various applications \cite{PR05, DP13, DP17}. We will consider both IMEX Runge-Kutta schemes such as IMEX-ARS scheme \cite{ARS97} and IMEX multistep schemes such as IMEX-BDF scheme \cite{HR07}. These propagation schemes together with the RNNs will make up our ``densely connected recurrent neural network" (DC-RNN).


Our third main contribution is to propose physics-based regularization to the loss function in \eqref{eq:intro} to improve optimization efficiency and avoid over-fitting. First, a physically correct model is usually described with a small number of mathematical operators in \eqref{eq:rnnterms2}, while an over-fitting model would have a large number of operators for a better fitting capacity. Thus, inspired by the lasso approaches in \cite{Tibshirani1996RegressionSA,BSPJ15,8573778}, we propose sparse regularization to avoid over-fitting and remove undesirable features in the governing equation, e.g., adding a $L^1$-norm penalty term to the coefficients in \eqref{eq:rnnterms2}. Second, a micro-macro decomposition of kinetic equations \cite{LM:08} is applied to transfer a challenging recovery problem with a single PDE to an easier recovery problem with a coupled PDE system, enforcing our recovery results to be more physically meaningful. Furthermore, the microscopic part, denoted as $g$ (see Section~\ref{ModelEquation}), satisfies
\begin{equation}\label{eq:intro4}
\langle g\rangle:=\int_{[-1,1]}g(v,x,t)dv=0,
\end{equation}
which will be used as a constraint of our recovery. Finally, in most cases, kinetic equations have spatial-dependent coefficients, which motivates us to design spatial-dependent parameters $\bd{\theta}(x)$ in \eqref{eq:rnnterms2} and the regularity in terms of $x$ can also be considered as a regularization penalty.

To summarize, the main highlights of our learning algorithm are as follows:

\begin{itemize}
\item DC-RNN built for transport (local) and  collision (nonlocal) operators typically involved in kinetic equations.
\item Multiscale-aware RNN structures and learning rates for the recovery of time-dependent PDEs.
\item Novel optimization loss function inspired by high-order IMEX for multiscale equations.
\item Physics-aware loss function and regularization specialized for kinetic equations.
\item Efficient arithmetic and memory cost.
\end{itemize}

We structure this manuscript as follows. In Section 2, an exemplary PDE for our learning problem is introduced to motivate our algorithm. In Section 3, we mathematically formulate an ansatz that we will use to fit data to PDEs. In Section 4, our physics-aware loss function is introduced to learn PDEs from data. In Section 5, we will carry out several numerical experiments to test our algorithm. Finally, concluding remarks are made in Section 6. 


\section{Model Equation: the Linear Transport Kinetic Equation}
\label{ModelEquation}

We now present a model equation, the linear transport equation, to motivate our learning algorithm. The linear transport equation is a prototype kinetic equation describing particles such as neutrons or photons interacting with a background medium \cite{Chandrasekhar, Davison}. This equation highlights some of the challenging aspects that an efficient learning algorithm should account for. That is, our model equation will allow us to understand the hypothesis space (the set of functions describing kinetic equations) better. This will lead us to devise ways to capture multiple scales, nonlocal operators, and regularity conditions. In addition, we will be able to discern appropriate numerical techniques needed to carry out our learning algorithm.

\comment{
We now present a model equation, the linear transport equation, to motivate our learning algorithm. The linear transport equation is a prototype kinetic equation describing particles such as neutrons or photons interacting with a background medium \cite{Chandrasekhar, Davison}. This equation highlights some of the challenging aspects that an efficient learning algorithm should account for: multiple scales in both time and space, a nonlocal integral operator, and a meaningful macroscopic limit.}

In the simple 1D case, the linear transport equation reads
\begin{equation}\label{eq:linear_transport}
\partial_{t}f + \dfrac{1}{\veps}v\partial_{x}f=\dfrac{\sigma^{S}}{\veps^{2}}(\langle f\rangle-f)-\sigma^{A}f +G,
\end{equation}
where $f=f(t,x,v)$ is the probability density function of time $t\geq 0$, position $x\in \Omega \subset \mathbb{R}$, and velocity $v\in [-1,1]$; $\langle\cdot\rangle:=\dfrac{1}{2}\int_{-1}^1 \cdot \rd{v}$ is a projection or collision operator; $\sigma^{S}(x)$ and $\sigma^{A}(x)$ are the scattering and absorption coefficients; and $G(x)$ is a given source. Finally, $\veps$ is a dimensionless parameter indicating the strength of the scattering. Indeed, when $\veps\sim \mathcal{O}(1)$, the equation (\ref{eq:linear_transport}) is in the fully kinetic regime (all operators balance); when $\veps \rightarrow 0$, the scattering is so strong that (\ref{eq:linear_transport}) approaches a diffusion limit. To see this, consider the so-called micro-macro decomposition of $f$:
\begin{equation}\label{eq:decomp}
f=\rho+\veps g, \quad \rho:=\langle f \rangle,
\end{equation}
where $\rho$ is the macro part (density) of the solution, and $g$ is the micro part. A crucial condition we use is
\begin{equation}\label{eq:gnull}
\langle g \rangle = 0.
\end{equation}
Equation \eqref{eq:gnull} is the conservation condition and will be numerically indispensable since it allows us to impose exact conditions satisfied by kinetic equations. Substituting (\ref{eq:decomp}) into (\ref{eq:linear_transport}), one can derive the following coupled system for $\rho$ and $g$, equivalent to \eqref{eq:linear_transport}:
\begin{align}
\partial_{t}\rho&=-\partial_{x}\langle vg\rangle-\sigma^{A}\rho + G, \label{eq:rho}\\
\partial_{t}g&=-\dfrac{1}{\veps}\left(\mathcal{I}-\langle \ \rangle\right)(v\partial_{x}g)-\dfrac{1}{\veps^{2}}v\partial_{x}\rho-\dfrac{\sigma^{S}}{\veps^{2}}g-\sigma^{A}g, \label{eq:GGG}
\end{align}
where $\mathcal{I}$ denotes the identity operator.

In (\ref{eq:GGG}), if $\veps\rightarrow 0$, one obtains
\begin{equation}
g= -\frac{1}{\sigma^S}v\partial_x\rho+\mathcal{O}(\veps),
\end{equation}
which, when substituted into (\ref{eq:rho}), yields
\begin{equation}\label{eq:diffusionlimit}
\partial_{t}\rho =\partial_{x}\left(\dfrac{1}{3\sigma^{S}}\partial_{x}\rho\right)-\sigma^{A}\rho +G +\mathcal{O}(\veps). 
\end{equation}
So $\rho$ follows the dynamics of a diffusion equation. We now go through a few things that we can learn from the linear transport equation following the notations used in Section 1.

$\bullet$ \textbf{Involved Basic Mathematical Operators: Identity, Advection, and Projection.} Notice that each of Equations \eqref{eq:linear_transport}, \eqref{eq:rho}, and \eqref{eq:GGG} can be recovered from the ansatz:

\begin{equation}\label{eq:smallansatz}
\begin{split}
\partial_{t}u&=a_{1}\mathcal{A}_{1}(g)+a_{2}\mathcal{A}_{2}(g)+a_{3}\mathcal{A}_{3}(g)+\sum_{i=1}^{3}\sum_{j=1}^{3}a_{i,j}\mathcal{A}_{i}\circ\mathcal{A}_{j}(g)\\
&+b_{1}\mathcal{A}_{1}(\rho)+b_{2}\mathcal{A}_{2}(\rho)+b_{3}\mathcal{A}_{3}(\rho)+\sum_{i=1}^{3}\sum_{j=1}^{3}b_{i,j}\mathcal{A}_{i}\circ\mathcal{A}_{j}(\rho)+B
\end{split}
\end{equation}
for $u=f$, $g$, or $\rho$. 
For example, the equation for $g(v,x,t)$ can be recovered provided
\begin{equation}\label{eq:basicops}
\mathcal{A}_{1}=\mathcal{I},\hspace{1cm}\mathcal{A}_{2}=v\partial_{x},\hspace{1cm}\mathcal{A}_{3}=\langle\cdot\rangle,
\end{equation}
with coefficients
\begin{equation}
\begin{split}
a_{1}=\dfrac{\sigma^{S}(x)}{\veps^{2}}-\sigma^{A}(x),\hspace{1cm}a_{2}=-\dfrac{1}{\veps},&\hspace{1cm}a_{3,2}=\dfrac{1}{\veps},\hspace{1cm}b_{2}=-\dfrac{1}{\veps^2}\\
a_{3}=0,\hspace{1cm}a_{i\neq 3,j\neq 2}=0,&\hspace{1cm}b_{i,j}=0,\hspace{1cm}B=0.\\
\end{split}
\end{equation}

Thus, at the very minimum, our hypothesis space in Equation \eqref{eq:intro2} should involve operators in  \eqref{eq:basicops}. We expect to see these operators for general kinetic equations. Potentially one can also have cubic or higher order nonlinearities in our hypothesis space. Therefore, we want to generalize Equation \eqref{eq:smallansatz} to involve greater number of compositions.

$\bullet$ \textbf{Functions of $x$.} Equations \eqref{eq:linear_transport}, \eqref{eq:rho}, and \eqref{eq:GGG} involve the functions $\sigma^{A}(x)$, $\sigma^{S}(x)$, and $G(x)$. Therefore the coefficients $\{a_{i},a_{i,j},b_{i},b_{i,j},B\}$ should be allowed to depend on $x$.

$\bullet$ \textbf{Scale Disparity.} If we want to determine the correct order of each term, then we need to make an asymptotic expansion:
\begin{equation}\label{eq:scaledisparity}
\begin{split}
a_{i} &= a_{i}^{0}+\dfrac{1}{\veps}a_{i}^{1}+\dfrac{1}{\veps^{2}}a_{i}^{2},\hspace{1cm}a_{i,j} = a_{i,j}^{0}+\dfrac{1}{\veps}a_{i,j}^{1}+\dfrac{1}{\veps^{2}}a_{i,j}^{2},\\
b_{i} &= b_{i}^{0}+\dfrac{1}{\veps}b_{i}^{1}+\dfrac{1}{\veps^{2}}b_{i}^{2},\hspace{1cm} b_{i,j} = b_{i,j}^{0}+\dfrac{1}{\veps}b_{i,j}^{1}+\dfrac{1}{\veps^{2}}b_{i,j}^{2},\\
B &= B^{0}+\dfrac{1}{\veps}B^{1}+\dfrac{1}{\veps^{2}}B^{2},\\
\end{split}
\end{equation}
where it is understood that the upper index labels the order of the scale. The multiscale phenomenon here is the main motivation of the multiscale model in Equation \eqref{eq:intro2}.

$\bullet$ \textbf{Exact Conditions.} Typically, adding regularization to machine learning problems can vastly improve the outcome of the prediction. There is one obvious constraint for our target kinetic equation: Equation \eqref{eq:gnull}. An added feature about this condition is that it is independent of $\veps$ and thus helpful for modeling dynamics between the small and large scale limits. 

$\bullet$ \textbf{Sparsity.} The large number of basis terms in our hypothesis space means that we might have overfitting issues. Thus, the following sparsity regularization term could be considered:
\begin{equation}
\begin{split}
\sum_{n}&\left((\sum_{i}||a_{i}^{n}||_{L^{1}}+||b_{i}^{n}||_{L^{1}}+||B^{n}||_{L^{1}})
+(\sum_{i,j}||a_{i,j}^{n}||_{L^{1}}+||b_{i,j}^{n}||_{L^{1}})\right),\\
\end{split}
\end{equation}
which can be enforced by adding the following regularization term to our loss function in Equation \eqref{eq:intro}:
\begin{equation}\label{eq:sparseiff}
||\bd{\theta}(x)||_{L^{1}},
\end{equation}
since $\bd{\theta}(x)$ is the actual parameters to be optimized in our model in Equation \eqref{eq:intro2}.

$\bullet$ \textbf{Numerical Techniques.} Finally, we need to consider numerical methods for arriving at the correct set of trainable parameters. For this reason, we will use the IMEX schemes for propagating $u(t_{n})\rightarrow u(t_{n+1})$, where the time rate of change is given by an ansatz like Equation \eqref{eq:smallansatz}. We will use gradient descent, specifically the Adam algorithm, to update trainable parameters.

\vspace{0.25cm}

In sum, the discussion above illustrates the motivation of our optimization problem, model design, and regularization terms introduced in Section 1.

\section{Formulating an Ansatz to Fit Data to Kinetic Equations}
\label{thirdsection}
In this section, we will construct an ansatz capable of representing Equations \eqref{eq:rho}, \eqref{eq:GGG}. For simplicity, let us focus on the case when the spatial variable $x$ is one-dimensional. It is easy to generalize to high-dimensional cases.  We start by introducing notations which will be used throughout the paper.

\textbf{Notation.} The functions involved in \eqref{eq:rho} or \eqref{eq:GGG} are multidimensional, e.g. $\rho=\rho(x,t)$ and $g=g(v,x,t)$. The values of $\rho$ and $g$ will be defined on a mesh $(x_{j},t_{k})$ and $(v_{i},x_{j},t_{k})$ for $i\in 1,\cdots N_{v}$, $j\in 1,\cdots N_{x}$, and $k\in 1,\cdots N_{t}$. To further simplify the notation, we will use $u:=u_{i,j}$ to denote $u$ as a scalar function evaluated at the $(i,j)$-th position corresponding to $(v_{i},x_{j})$. The upper index $n$ in $u^{n}:=u(\cdot,t_{n})$ will correspond to time with $u^{1:N_t}:=\{u(\cdot,t_{i});\,\,\text{for}\,\, i\in 1\cdots N_t\}$ denoting $u$ evaluated at a time sequence. Matrices will be written with capital letters while operators applied to the data will mainly be written using script letters.

\subsection{Operator Evaluation}\label{sec:OE}

We will describe the evaluation of commonly used operators in Equations \eqref{eq:rho}, \eqref{eq:GGG}, and other kinetic equations below.

\emph{1) Identity operator.}
The identity operator is defined by 
\begin{equation}
\mathcal{I}(u):=u.
\end{equation}
The evaluation of $\mathcal{I}(u)$ at the point $(v_{i},x_{j})$ simply follows 
$\mathcal{I}(u)_{i,j}:=u_{i,j}$.

\emph{2) Pseudo-upwind for the advection operator.}
We define the advection operator acting on $u$ as the dot-product:
\begin{equation}\label{eq:theadvection}
{v}\cdot\nabla_{{x}}u,
\end{equation}
where ${v}$ is a velocity distribution. We note that many stable schemes use an upwind stencil for the advection operator. The first-order upwind stencil gives:
\begin{equation*}
\partial_{x}u_{i,j}^{+}=\dfrac{u_{i,j+1}-u_{i,j}}{\Delta x}\hspace{0.5cm}\text{for}\,\,v>0,\\
\end{equation*}
\begin{equation*}
\partial_{x}u_{i,j}^{-}=\dfrac{u_{i,j}-u_{i,j-1}}{\Delta x}\hspace{0.5cm}\text{for}\,\,v<0,\\
\end{equation*}
and
\begin{equation*}
v\partial_{x}u_{i,j}=v_{-}\partial_{x}u_{i,j}^{+}+v_{+}\partial_{x}u_{i,j}^{-},
\end{equation*}
which is the evaluation of the advection operator in \eqref{eq:theadvection} at the point $(v_{i},x_{j})$ in the one-dimensional case. This stencil is suitable for a first-order-in-time IMEX-scheme. For higher-order IMEX schemes, one should use higher-order stencils.

\emph{3) Projection operators.}
We define the projection operator as an integral with respect to the variable $v$ of a function $u(v,x)$. In one-dimension, we have:
\begin{equation}\label{eq:theprojection}
\langle u\rangle:=\dfrac{1}{2}\int_{-1}^{1}u(v,x) \,\,dv,
\end{equation}
which can be discretized as a finite sum using Gauss quadrature:
\begin{equation*}
\langle u\rangle_j\approx\dfrac{1}{2}\sum_{i=1}^{N_{v}}u_{i,j}\cdot w_{i}
\end{equation*}
with quadrature weights $\{w_{i}\}$. Note that the data corresponding to $u$ is represented by a two-index tensor $u_{i,j}$ with $i,j$ corresponding to the values $v_{i}$ and $x_{j}$, respectively. The above quadrature maps $u$ to a one-index tensor $\langle u\rangle_j$. To make dimensions consistent, we extend this to a two-index tensor by $\langle u \rangle_{i,j}:=\langle u\rangle_{j}$ for each $i$.

\emph{4) Other differential operators.}
Higher-order differential operators such as the Laplacian will be computed by using central difference formulas.

\subsection{Ansatz for Fitting PDEs to Data}

In this section, we will form an ansatz that will be used to fit a PDE to data, i.e., identifying the governing PDE to which the observed data is a discrete solution. We will consider the following two typical examples for simplicity. The generalization to other cases is simple.

\textbf{Scalar equation ansatz.}
Let us consider a first-order in time PDE, then the equation ansatz is built as
\begin{equation}\label{eq:PTexpansionMain}
\partial_{t}u=\mathcal{F}(u)
\end{equation}
with $F$ split into $M$ multiscale components following our main model in \eqref{eq:intro2}:
\begin{equation*}
\mathcal{F}:= \mathcal{F}^{0}(u)+\dfrac{1}{\veps_{pred}}\mathcal{F}^{1}(u)+\dfrac{1}{\veps_{pred}^{2}}\mathcal{F}^{2}(u)+\cdots+\dfrac{1}{\veps_{pred}^{M}}\mathcal{F}^{M}(u).
\end{equation*}
The integer $M$ depends on the number of multiscale components for the problem being considered. If one only expects one fast scale and one slow scale component, $M$ is set to $M=1$. For slow, medium, and fast scales, $M$ is set to $M=2$, etc. $\veps_{pred}$ is a learnable scaling number defined in \eqref{eqn:eps} restricted to $0<\veps_{pred}\leq 1$ but not necessarily equal to $\veps$. The operators $\mathcal{F}^{0}, \mathcal{F}^{1}, \mathcal{F}^{2}, \cdots$  will be differential operators acting on $u$ and constructed as in \eqref{eq:rnnterms2}. The construction detail will be provided in the next section.

\textbf{Two-component vector equation ansatz.}
For vectorized equations, we build an ansatz for each component individually as
\begin{equation}\label{eq:vectorizedansatz}
\begin{split}
\partial_{t}g=\mathcal{F}_{1}(g,\rho)=\sum_{m=0}^{M}\dfrac{1}{\veps_{pred}^{m}}(\mathcal{F}_{1,1}^{m}(g)+\mathcal{F}_{1,2}^{m}(\rho))\\
\partial_{t}\rho=\mathcal{F}_{2}(g,\rho)=\sum_{m=0}^{M}\dfrac{1}{\veps_{pred}^{m}}(\mathcal{F}_{2,1}^{m}(g)+\mathcal{F}_{2,2}^{m}(\rho)).\\
\end{split}
\end{equation}
The $\mathcal{F}_{q,p}^{m}$ are generally different operators for each $q$,$p$, and $m$ following the construction in \eqref{eq:rnnterms2}. Each $\mathcal{F}_{q,p}^{m}$ has an individual set of network parameters. $\veps_{pred}$ is a learnable scaling number as in the previous example. For the remainder of the manuscript, $\mathcal{F}_{1}$ will denote the right hand side of the $g$-equation. $\mathcal{F}_{2}$ will denote the right hand side of the $\rho$-equation. The construction of $\mathcal{F}_{q,p}^{m}$ using an RNN structure will be presented in detail in the next section.

\remark Equation \eqref{eq:vectorizedansatz} is our chosen ansatz. There are many alternative ways to construct an ansatz. For example, we only consider the linear combination of $\mathcal{F}_{q,p}^{m}$ and it is also possible to explore the products of $\mathcal{F}_{q,p}^{m}$. 

\subsection{Building a Dictionary Using RNNs}\label{dictionarybuilding}

We construct the operators $\mathcal{F}^{m}$ in Equation \eqref{eq:PTexpansionMain} for the single-component case. For the multicomponent case as in \eqref{eq:vectorizedansatz}, we construct $\mathcal{F}_{q,p}^{m}$ in the same manner. The only difference is that each $F_{q,p}^{m}$ will have a different set of parameters depending on $(q,p)$. We will omit the $(q,p)$ index for clarity. To begin with, we will need to supply the RNN with a few basic mathematical operators as mentioned in the introduction. In particular, we consider operators $\mathcal{A}_{1}(u)=\mathcal{I}(u)$, $\mathcal{A}_{2}(u)=v\cdot \nabla u$, and $\mathcal{A}_{3}(u)=\langle u\rangle$ discussed in Section \ref{sec:OE}. It is potentially better to include more basic mathematical operators such as $\mathcal{A}_{4}(u)=\nabla u$, $\mathcal{A}_{5}(u)=g^{2}\cdot u$, $\mathcal{A}_{6}(u)=\exp\left(-(u)^{2}\right)$, etc. 

Next, a symbolic RNN will be introduced to generate a complicated operator $\mathcal{F}^{m}$ using basic mathematical operators. Given basic mathematical operators $\{\mathcal{A}_{1},\cdots,\mathcal{A}_{n}\}$, we build a $k$-layer RNN for each $m=0,1,2,\cdots,$ by successively applying a weight matrix $W_{2,n}\in\mathbb{R}^{2\times n}$ to the operator vector $[\mathcal{A}_{1},\cdots,\mathcal{A}_{n}]^{T}$ and then adding a bias vector $B_{2,1}=[b_1,b_2]^T\in\mathbb{R}^2$ times $\mathcal{I}$:
\begin{equation}\label{eq:ope1}
\begin{split}
W_{2,n}[\mathcal{A}_{1},\cdots,\mathcal{A}_{n}]^{T}+B_{2,1}\mathcal{I}&:=
\left[\begin{matrix}
w_{1,1}\mathcal{A}_{1}+w_{1,2}\mathcal{A}_{2}+\cdots+w_{1,n}\mathcal{A}_{n}+b_{1}\mathcal{I}\\
w_{2,1}\mathcal{A}_{1}+w_{2,2}\mathcal{A}_{2}+\cdots+w_{2,n}\mathcal{A}_{n}+b_{2}\mathcal{I}\\
\end{matrix}
\right]\\
&:=
\left[\begin{matrix}
\mathcal{C}_{1}\\
\mathcal{C}_{2}\\
\end{matrix}
\right].
\end{split}
\end{equation}
Because $\mathcal{C}_{1}$ and $\mathcal{C}_{2}$ are operators, they can be applied to generate a more expressive formulation with a special ``composition" denoted as $\odot$ defined below:
\comment{
\begin{equation*}
\begin{split}
\mathcal{C}_{1}\odot\mathcal{C}_{2}&:=w_{1,1}w_{2,1}\mathcal{A}_{1}\circ\mathcal{A}_{1}+w_{1,1}w_{2,2}\mathcal{A}_{1}\circ\mathcal{A}_{2}+\cdots+ w_{1,1}w_{2,n}\mathcal{A}_{1}\circ\mathcal{A}_{n}+\cdots\\
&+w_{1,n}w_{2,1}\mathcal{A}_{n}\circ\mathcal{A}_{1}+w_{1,n}w_{2,2}\mathcal{A}_{n}\circ\mathcal{A}_{2}+\cdots+ w_{1,n}w_{2,n}\mathcal{A}_{n}\circ\mathcal{A}_{n}\\
&+(w_{1,1}b_{2}+w_{2,1}b_{1})\mathcal{A}_{1}+(w_{1,2}b_{2}+w_{2,2}b_{1})\mathcal{A}_{2}+\cdots+(w_{1,n}b_{2}+w_{2,n}b_{1})\mathcal{A}_{n}
\end{split}
\end{equation*}
}

\begin{equation}\label{eq:ope2}
\begin{split}
\mathcal{C}_{1}\odot\mathcal{C}_{2}&:=w_{1,1}w_{2,1}\mathcal{A}_{1}\circ\mathcal{A}_{1}+\cdots+ w_{1,1}w_{2,n}\mathcal{A}_{1}\circ\mathcal{A}_{n}+\cdots\\
&+w_{1,n}w_{2,1}\mathcal{A}_{n}\circ\mathcal{A}_{1}+\cdots+ w_{1,n}w_{2,n}\mathcal{A}_{n}\circ\mathcal{A}_{n}\\
&+(w_{1,1}b_{2}+w_{2,1}b_{1})\mathcal{A}_{1}+\cdots+(w_{1,n}b_{2}+w_{2,n}b_{1})\mathcal{A}_{n},
\end{split}
\end{equation}
where $\circ$ denotes the standard composition.


Now we define $\mathcal{F}^{m}$:

\begin{equation}
\begin{split}
\bd{\xi}^{(1)}:&=W_{2,n}^{1,m}[\mathcal{A}_{1},\mathcal{A}_{2},\cdots,\mathcal{A}_{n}]^{T}+B_{2,1}^{1,m}\mathcal{I}\\
\mathcal{B}_{1}:&=\mathcal{C}_{1}^{(1)}\odot\mathcal{C}_{2}^{(1)}\\
\bd{\xi}^{(2)}:&=W_{2,n+1}^{2,m}[\mathcal{A}_{1},\mathcal{A}_{2},\cdots,\mathcal{A}_{n},\mathcal{B}_{1}]^{T}+B_{2,1}^{2,m}\mathcal{I}\\
\mathcal{B}_{2}:&=\mathcal{C}_{1}^{(2)}\odot\mathcal{C}_{2}^{(2)}\\
\bd{\xi}^{(3)}:&=W_{2,n+2}^{3,m}[\mathcal{A}_{1},\mathcal{A}_{2},\cdots,\mathcal{A}_{n},\mathcal{B}_{1},\mathcal{B}_{2}]^{T}+B_{2,1}^{3,m}\mathcal{I}\\
\mathcal{B}_{3}:&=\mathcal{C}_{1}^{(3)}\odot\mathcal{C}_{2}^{(3)}\\
&\vdots\\
\bd{\xi}^{(K)}:&=W_{2,n+K-1}^{K,m}[\mathcal{A}_{1},\mathcal{A}_{2},\cdots,\mathcal{A}_{n},\mathcal{B}_{1},\cdots,\mathcal{B}_{K-1}]^{T}+B_{2,1}^{K,m}\mathcal{I}\label{eq:eq210}\\
\mathcal{B}_{K}:&=\mathcal{C}_{1}^{(K)}\odot\mathcal{C}_{2}^{(K)}\\
\mathcal{F}^{m}:&=W_{1,n+K}^{K+1,m}[\mathcal{A}_{1},\mathcal{A}_{2},\cdots,\mathcal{A}_{n},\mathcal{B}_{1},\cdots,\mathcal{B}_{K}]^{T},
\end{split}
\end{equation}
where the weight matrices are given by:
\begin{equation}\label{eq:finalcoefficients}
W_{2,n+k-1}^{k,m}:=\left[\begin{matrix}
w_{1,1}^{k,m} & w_{1,2}^{k,m} & \cdots & w_{1,n+k-1}^{k,m}\\
w_{2,1}^{k,m} & w_{2,2}^{k,m} & \cdots & w_{2,n+k-1}^{k,m}\\
\end{matrix}\right]
\end{equation}
for $k=1,\cdots,K$ and,
\begin{equation*}
W_{1,n+K}^{K+1,m}:=\left[\begin{matrix}
w_{1}^{K+1,m} & w_{2}^{K+1,m} & \cdots & w_{n+K}^{K+1,m}\\
\end{matrix}\right]
\end{equation*}
with each $w_{i,j}^{k,m}\in \mathbb{R}$. The biases are given by:
\begin{equation*}
B_{2,1}^{k,m}:=\left[\begin{matrix}
b_{1}^{k,m}\\
b_{2}^{k,m}\\
\end{matrix}\right]
\end{equation*}
with $b_{j}^{k,m}\in \mathbb{R}$.

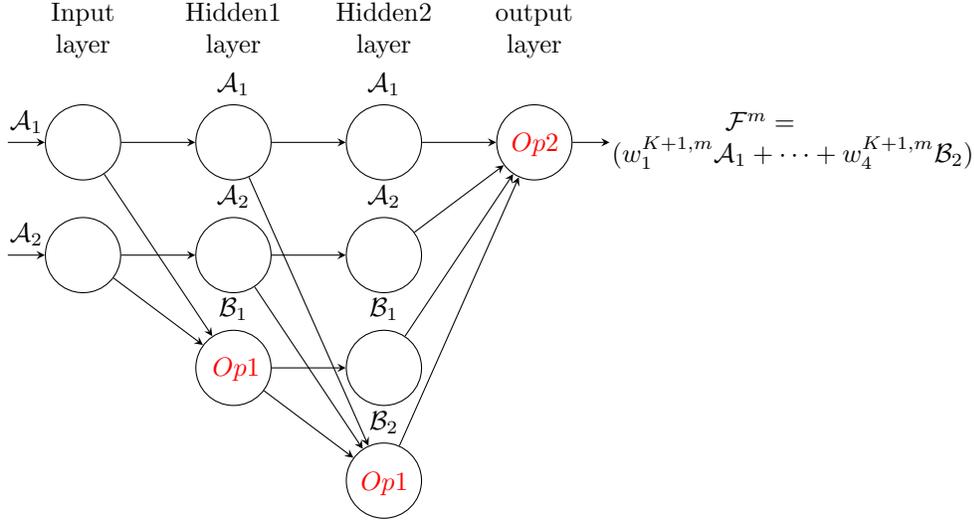
\begin{figure}\label{NNfig}
\tikzset{%
  every neuron/.style={
    circle,
    draw,
    minimum size=1cm
  },
  neuron missing/.style={
    draw=none, 
    scale=4,
    text height=0.033cm,
    execute at begin node=\color{black}$\vdots$
  },
}

\begin{tikzpicture}[x=1.0cm, y=1.0cm, >=stealth]

\foreach \m/\l [count=\y] in {1,2}
  \node [every neuron/.try, neuron \m/.try] (input-\m) at (0,2.5-\y*1.5) {};

\foreach \m [count=\y] in {1,2,3}
  \node [every neuron/.try, neuron \m/.try ] (hidden1-\m) at (2,2.5-\y*1.5) {};

\foreach \m [count=\y] in {1,2,3,4}
  \node [every neuron/.try, neuron \m/.try ] (hidden2-\m) at (4,2.5-\y*1.5) {};

\foreach \m [count=\y] in {1}
  \node [every neuron/.try, neuron \m/.try ] (output-\m) at (6,2.5-\y*1.5) {};

\put(49, -60){$\color{red}Op1$}
\put(105, -103){$\color{red}Op1$}
\put(162, 25){$\color{red}Op2$}
\put(243, 33){$\mathcal{F}^{m}=$};
\put(200, 21){$(w_{1}^{K+1,m}\mathcal{A}_{1}+\cdots+w_{4}^{K+1,m}\mathcal{B}_{2})$};

\foreach \l [count=\i] in {1,2}
  \draw [<-] (input-\i) -- ++(-1,0)
    node [above, midway] {$\mathcal{A}_\l$};

\foreach \l [count=\i] in {1,2}
 \node [above] at (hidden1-\i.north) {$\mathcal{A}_\l$};
\foreach \l [count=\i from 3] in {1}
  \draw [->] (hidden1-\i)
    node [above] at (hidden1-\i.north) {$\mathcal{B}_{\l}$};

\foreach \l [count=\i] in {1,2}
  \draw [->] (hidden2-\i)
    node [above] at (hidden2-\i.north) {$\mathcal{A}_\l$};
\foreach \l [count=\i from 3] in {1,2}
  \draw [->] (hidden2-\i)
    node [above] at (hidden2-\i.north){$\mathcal{B}_{\l}$};

\foreach \l [count=\i] in {1}
  \draw [->] (output-\i);

\foreach \i in {1,2}
  \foreach \j in {3}
    \draw [->] (input-\i) -- (hidden1-\i);

\foreach \i in {1,2}
  \foreach \j in {3}
    \draw [->] (input-\i) -- (hidden1-\j);

\foreach \i in {1,...,3}
  \foreach \j in {4}
    \draw [->] (hidden1-\i) -- (hidden2-\j);

\foreach \i in {1,2,3}
    \draw [->] (hidden1-\i) -- (hidden2-\i);

\foreach \i in {1,...,4}
\foreach \j in {1}
    \draw [->] (hidden2-\i) -- (output-\j);

\foreach \l [count=\i] in {1} 
  \draw [->] (output-\i) -- ++(1,0);

\foreach \l [count=\x from 0] in {Input, Hidden1, Hidden2,output}
  \node [align=center, above] at (\x*2,2) {\l \\ layer};

\end{tikzpicture}
\caption{Example RNN with $K=2$ hidden layers. The output $\mathcal{F}^{m}$ makes up the order $\veps^{m}$ part of the right hand side of the PDE. $Op1$ is the mathematical operation in Equation \eqref{eq:ope1} and \ref{eq:ope2} which takes linear combination plus bias of the previous layer and performing a composition. $Op2$ is the operation of forming a linear combination of the previous layer (the last line of Equation \eqref{eq:eq210}).}
\end{figure}

The operator built by the recursive compositions in \eqref{eq:eq210} is a symbolic RNN operator, the evaluation of which on a given function follows in the basic evaluation rules introduced in Section \ref{sec:OE}. will have to be evaluated at both data sets $\{g(v,x,t_{i})\}$ and $\{\rho(x,t_{i})\}$, since our model problem depends on both $g$ and $\rho$. A diagrammatic representation of this RNN is shown in Figure \ref{NNfig}.

\remark{We have adopted the recursive framework introduced in \cite{BinD18} to build our RNN. The main difference between the RNN in this manuscript and the RNN in \cite{BinD18} is that our RNN can learn nonlocal and multiscale operators. Other RNN frameworks may also be good alternatives. Optimizing the RNN framework is not a focus in this paper. }


\remark The weights and biases can be trainable
space-dependent functions such that our algorithm can learn more space-dependent operators, e.g., let
\begin{equation}
w_{i,j}^{k,m}(x):\mathbb{R}\rightarrow\mathbb{R},\,\,\,\text{and}\,\,\,b_{j}^{k,m}(x):\mathbb{R}\rightarrow\mathbb{R}.
\end{equation}
In more particular, one can also replace these weights and biases with neural networks in the spatial variable $x$ at the cost of using more parameters. We will let the reader explore these possibilities but, we will also present a yet different alternative to treating space-dependent weights and biases in the next section.

\subsection{An Example when $K=1$}
Using $K=1$ and two basic mathematical operators $\mathcal{A}_{1}$ and $\mathcal{A}_{2}$ in \eqref{eq:eq210}, for the PDE model in \eqref{eq:PTexpansionMain}, we produce an RNN as a scalar PDE ansatz of the form:
\begin{equation}
\begin{split}
\partial_{t}u=&\sum_{m=0}^{M}\dfrac{1}{\veps_{pred}^{m}}\left(w_{1}^{2,m}\mathcal{A}_{1}+w_{2}^{2,m}\mathcal{A}_{2}+w_{3}^{2,m}\left[(w_{1,1}^{1,m}\mathcal{A}_{1}+w_{1,2}^{1,m}\mathcal{A}_{2}\right.\right.\\
+&\left.\left. b_{1}^{1,m}I.d.)\circ (w_{2,1}^{1,m}\mathcal{A}_{1}+w_{2,2}^{1,m}\mathcal{A}_{2}+b_{2}^{1,m}I.d.)\right]\right)(u)\\
=&\sum_{m=0}^{M}\dfrac{1}{\veps_{pred}^{m}}\left[w_{3}^{2,m}b_{1}^{1,m}b_{2}^{1,m}u+(w_{1}^{2,m}+w_{3}^{2,m}(w_{1,1}^{1,m}b_{2}^{1,m}+b_{1}^{1,m}w_{2,1}^{1,m}))\mathcal{A}_{1}(u)\right.\\
+&\left.(w_{2}^{2,m}+w_{3}^{2,m}(w_{1,2}^{1,m}b_{2}^{1,m}+b_{1}^{1,m}w_{2,1}^{1,m}))\mathcal{A}_{2}(u)\right.\\
+&\left.w_{3}^{2,m}(w_{1,1}^{1,m}w_{2,1}^{1,m}\mathcal{A}_{1}\circ\mathcal{A}_{1}(u)+w_{1,1}^{1,m}w_{2,2}^{1,m}\mathcal{A}_{1}\circ\mathcal{A}_{2}(u))\right.\\
+&\left.w_{3}^{2,m}(w_{1,2}^{1,m}w_{2,1}^{1,m}\mathcal{A}_{2}\circ\mathcal{A}_{1}(u)+w_{1,2}^{1,m}w_{2,2}^{1,m}\mathcal{A}_{2}\circ\mathcal{A}_{2}(u))\right],
\end{split}
\end{equation}
with $\veps_{pred}$ given by Equation \eqref{eq:nonparallel} or \eqref{eq:parallel}. 

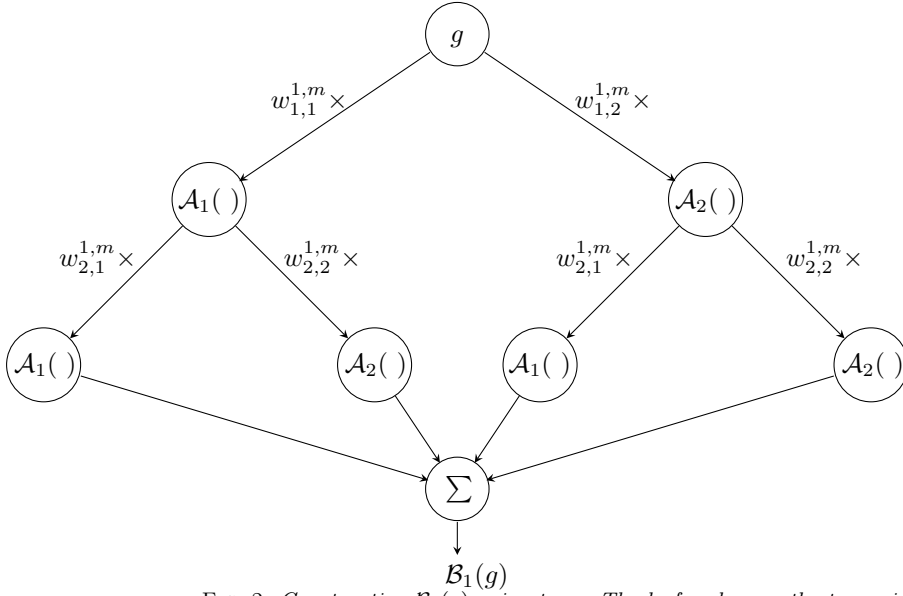
\begin{figure}\label{coeff_tree1}
\begin{tikzpicture}[x=1.1cm, y=1.1cm, >=stealth]
\draw[] (12,6.0) circle (12pt);
\put(373,185){$g$}
\draw [black, ->          ] (11.67,5.78) -- (9.364,4.221);
\put(305,160){$w_{1,1}^{1,m}\times$}
\put(270,122){$\mathcal{A}_{1}(\,\,)$}
\draw[] (9,4.0) circle (14pt);
\draw [black, ->          ] (12.33,5.78) -- (14.636,4.221);
\put(420,160){$w_{1,2}^{1,m}\times$}
\put(458,122){$\mathcal{A}_{2}(\,\,)$}
\draw[] (15,4.0) circle (14pt);

\draw [black, ->          ] (8.68,3.68) -- (7.32,2.32);
\draw[] (7,2.0) circle (14pt);
\put(225,100){$w_{2,1}^{1,m}\times$}
\put(208,60){$\mathcal{A}_{1}(\,\,)$}
\draw [black, ->          ] (9.32,3.68) -- (10.68,2.32);
\draw[] (11,2.0) circle (14pt);
\put(310,100){$w_{2,2}^{1,m}\times$}
\put(333,60){$\mathcal{A}_{2}(\,\,)$}
\draw [black, ->          ] (14.68,3.68) -- (13.32,2.32);
\draw[] (13,2.0) circle (14pt);
\put(413,100){$w_{2,1}^{1,m}\times$}
\put(395,60){$\mathcal{A}_{1}(\,\,)$}
\draw [black, ->          ] (15.32,3.68) -- (16.68,2.32);
\draw[] (17,2.0) circle (14pt);
\put(500,100){$w_{2,2}^{1,m}\times$}
\put(521,60){$\mathcal{A}_{2}(\,\,)$}
\draw [black, ->          ] (7.45,1.865) -- (11.65,0.605);
\draw [black, ->          ] (11.25,1.625) -- (11.79,0.815);
\draw [black, ->          ] (12.75,1.625) -- (12.21,0.815);
\draw [black, ->          ] (16.55,1.865) -- (12.35,0.605);
\draw[] (12,0.5) circle (12pt);
\put(371,13){$\sum$}
\draw [black, ->          ] (12.0,0.1) -- (12.0,-.3);
\put(371,-20){$\mathcal{B}_{1}(g)$}

\end{tikzpicture}
\caption{Constructing $\mathcal{B}_{1}(g)$ using trees. The leaf nodes are the terms in the third row.}
\end{figure}

\begin{figure}\label{coeff_tree2}
\begin{tikzpicture}[x=1.1cm, y=1.1cm, >=stealth]
\draw[] (14,6.0) circle (12pt);
\put(435,185){$g$}

\draw [black, ->          ] (13.64,5.82) -- (10.4,4.2);
\put(350,160){$w_{1,1}^{2,m}\times$}
\put(302,122){$\mathcal{A}_{1}(\,\,)$}
\draw[] (10,4.0) circle (14pt);

\draw [black, ->          ] (14,5.62) -- (14,4.44);
\put(440,160){$w_{1,2}^{2,m}\times$}
\put(427,122){$\mathcal{A}_{2}(\,\,)$}
\draw[] (14,4.0) circle (14pt);

\draw [black, ->          ] (14.36,5.82) -- (17.6,4.2);
\put(500,160){$w_{1,3}^{2,m}\times$}
\put(552,122){$\mathcal{B}_{1}(\,\,)$}
\draw[] (18,4.0) circle (14pt);

\draw [black, ->          ] (9.805,3.61) -- (8.695,1.39);
\draw[] (8.5,1.0) circle (14pt);
\put(262,80){$w_{2,1}^{2,m}\times$}
\put(255,28){$\mathcal{A}_{1}(\,\,)$}

\draw [black, ->          ] (10,3.55) -- (10,1.45);
\draw[] (10,1.0) circle (14pt);
\put(284,55){$w_{2,2}^{2,m}\times$}
\put(301,28){$\mathcal{A}_{2}(\,\,)$}

\draw [black, ->          ] (10.195,3.61) -- (11.305,1.39);
\draw[] (11.5,1.0) circle (14pt);
\put(338,80){$w_{2,3}^{2,m}\times$}
\put(349,28){$\mathcal{B}_{1}(\,\,)$}

\draw [black, ->          ] (13.73,3.64) -- (12.77,2.36);
\draw[] (12.5,2.0) circle (14pt);
\put(390,100){$w_{2,1}^{2,m}\times$}
\put(380,60){$\mathcal{A}_{1}(\,\,)$}

\draw [black, ->          ] (14,3.54) -- (14,2.45);
\draw[] (14,2.0) circle (14pt);
\put(410,80){$w_{2,2}^{2,m}\times$}
\put(427,60){$\mathcal{A}_{2}(\,\,)$}

\draw [black, ->          ] (14.27,3.64) -- (15.23,2.36);
\draw[] (15.5,2.0) circle (14pt);
\put(460,100){$w_{2,3}^{2,m}\times$}
\put(474,60){$\mathcal{B}_{1}(\,\,)$}

\draw [black, ->          ] (17.805,3.61) -- (16.695,1.39);
\draw[] (16.5,1.0) circle (14pt);
\put(511,80){$w_{2,1}^{k,m}\times$}
\put(505,28){$\mathcal{A}_{1}(\,\,)$}

\draw [black, ->          ] (18,3.55) -- (18,1.45);
\draw[] (18,1.0) circle (14pt);
\put(534,60){$w_{2,2}^{k,m}\times$}
\put(552,28){$\mathcal{A}_{2}(\,\,)$}

\draw [black, ->          ] (18.195,3.61) -- (19.305,1.39);
\draw[] (19.5,1.0) circle (14pt);
\put(587,80){$w_{2,3}^{k,m}\times$}
\put(600,28){$\mathcal{B}_{1}(\,\,)$}

\draw [black, ->          ] (11.5,0.54) -- (13.65,0.0756);
\draw [black, ->          ] (10,0.54) -- (13.64,0.0486);
\draw [black, ->          ] (8.5,0.54) -- (13.615,0.0378);

\draw [black, ->          ] (12.77,1.64) -- (13.775,0.3);
\draw [black, ->          ] (14,1.56) -- (14,0.38);
\draw [black, ->          ] (15.23,1.64) -- (14.225,0.3);

\draw [black, ->          ] (16.5,0.54) -- (14.35,0.0756);
\draw [black, ->          ] (18,0.54) -- (14.36,0.0486);
\draw [black, ->          ] (19.5,0.54) -- (14.385,0.0378);

\draw[] (14,0) circle (12pt);
\put(433,-4){$\sum$}
\draw [black, ->          ] (14.0,-0.39) -- (14.0,-.7);
\put(435,-33){$\mathcal{B}_{2}(g)$}
\end{tikzpicture}
\caption{Constructing $\mathcal{B}_{2}(g)$ using trees. This tree depends on tree $\mathcal{B}_{1}$.}
\end{figure}
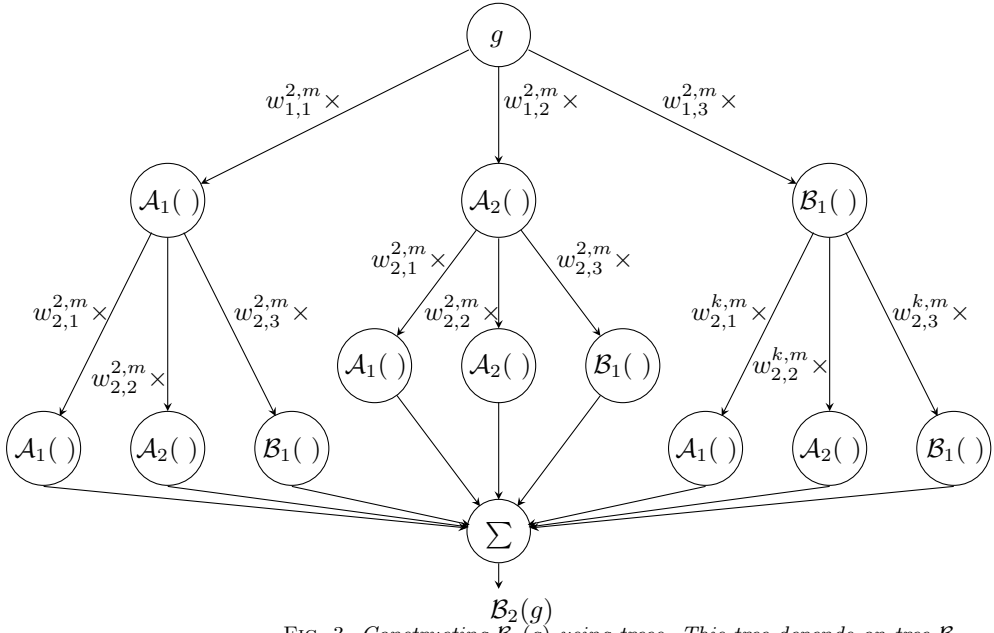

One way to calculate coefficients for the PDE is to construct trees for the operators $\mathcal{B}_{1}, \mathcal{B}_{2}, \cdots$. The benefit of using trees is that sub-trees following coefficients that are small in magnitude can be deleted. In figures \ref{coeff_tree1} and \ref{coeff_tree2}, we construct some example trees for $\mathcal{B}_{1}$ and $\mathcal{B}_{2}$ where we start with two initial operators $\mathcal{A}_{1}$ and $\mathcal{A}_{2}$ and omit the bias operator.

The weights and biases are determined by minimizing a loss function defined in the next section.

\section{Loss Functions for Learning PDEs}

To deduce the weights and biases for our PDE ansatz, we need to minimize a loss function. We begin by describing an unregularized loss function for learning PDEs from data.

\subsection{Unregularized Loss Function} Let us first focus on the case of a single scalar equation ansatz in \eqref{eq:PTexpansionMain}. 
We build an unregularized loss that will be a data-dependent function with the following abstract notation:
\begin{equation}\label{eq:unreg}
L(\bd{\theta})=\dfrac{1}{N_{t}-q}\sum_{n=1}^{N_t-q}||\mathcal{K}_{u}^{n}(D^{n,q};\bd{\theta})||_*
\end{equation}
where $\bd{\theta}$ denotes the set of all parameters in our RNN and 
\begin{equation}\label{eq:IMEXMap}
\mathcal{K}_{u}^{n}(D^{n,q};\bd{\theta})
\end{equation}
relates $q+1$-tuple data points:
\begin{equation*}
D^{n,q}:=\{u(x,t_{n}),u(x,t_{n+1}),\cdots,u(x,t_{n+q})\},\hspace{0.5cm}n=1,\cdots N_{t}-q.
\end{equation*}
The idea is that as $L\rightarrow 0$ with respect to a suitable norm $\|\cdot\|_*$, $\mathcal{F}$ approaches the correct PDE. Commonly used norms for loss minimization include $\ell_{1}$, $\ell_{2}$, and the Huber loss (see \cite{BinD18,DBLP:journals/corr/abs-1711-10561}).

To be precise, the relation of $D^{n,q}$ is specified by a time-stepping scheme, e.g., the IMEX scheme. However, to give the reader a greater understanding of $\mathcal{K}_{u}^{n}(D^{n,q};\bd{\theta})$, we will start with simpler schemes here. The symbolic RNN introduced in the previous section together with the IMEX schemes will make up our Densely Connected Recurrent Neural Network (DC-RNN).

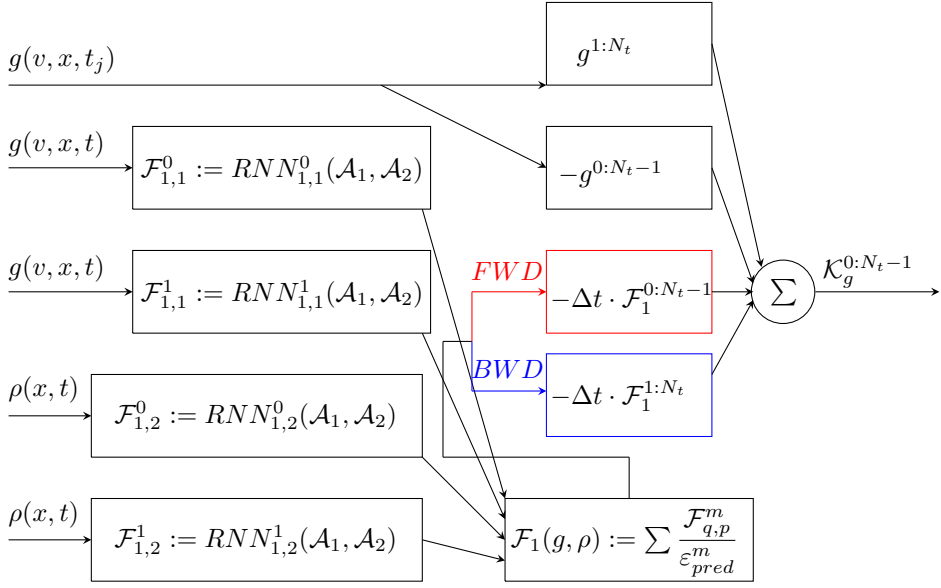
\begin{figure}\label{NNfigFwrd}
\begin{tikzpicture}[x=1.1cm, y=1.1cm, >=stealth]

\draw [black, ->          ] (0,4.5) -- (6.5,4.5);
\put(0,148){$g(v,x,t_{j})$}

\draw [black, ->          ] (4.5,4.5) -- (6.5,3.5);

\draw [black, ->          ] (0,3.5) -- (1.5,3.5);
\put(0,116){$g(v,x,t)$}
\put(50,106){$\mathcal{F}_{1,1}^{0}:=RNN_{1,1}^{0}(\mathcal{A}_{1}, \mathcal{A}_{2})$}
\draw[black] (1.5,3) -- (1.5,4) -- (5.1,4) -- (5.1,3) -- (1.5,3);

\draw [black, ->          ] (0,2) -- (1.5,2);
\put(0,69){$g(v,x,t)$}
\put(50,59){$\mathcal{F}_{1,1}^{1}:=RNN_{1,1}^{1}(\mathcal{A}_{1}, \mathcal{A}_{2})$}
\draw[black] (1.5,1.5) -- (1.5,2.5) -- (5.1,2.5) -- (5.1,1.5) -- (1.5,1.5);

\draw [black, ->          ] (0,0.5) -- (1,0.5);
\put(0,23){$\rho(x,t)$}
\put(40,13){$\mathcal{F}_{1,2}^{0}:=RNN_{1,2}^{0}(\mathcal{A}_{1}, \mathcal{A}_{2})$}
\draw[black] (1,0) -- (1,1) -- (5,1) -- (5,0) -- (1,0);

\draw [black, ->          ] (0,-1) -- (1,-1);
\put(0,-24){$\rho(x,t)$}
\put(40,-34){$\mathcal{F}_{1,2}^{1}:=RNN_{1,2}^{1}(\mathcal{A}_{1}, \mathcal{A}_{2})$}
\draw[black] (1,-1.5) -- (1,-0.5) -- (5,-0.5) -- (5,-1.5) -- (1,-1.5);

\put(190,-34){$\mathcal{F}_{1}(g,\rho):=\sum\dfrac{\mathcal{F}_{q,p}^{m}}{\veps_{pred}^{m}}$}
\draw (6,-1.5) -- (6,-0.5) -- (9,-0.5) -- (9,-1.5) -- (6,-1.5);

\draw [black, ->  ] (5,3) -- (6,-0.5);

\draw [black, ->  ] (5,1.5) -- (6,-0.75);

\draw [black, ->  ] (5,0) -- (6,-1);

\draw [black, ->  ] (5,-1.0) -- (6,-1.25);

\put(215,150){$g^{1:N_t}$}
\draw (6.5,4.5) -- (6.5,5.5) -- (8.5,5.5) -- (8.5,4.5) -- (6.5,4.5);

\put(208,104){$-g^{0:N_t-1}$}
\draw (6.5,3) -- (6.5,4) -- (8.5,4) -- (8.5,3) -- (6.5,3);

\put(205,20){$-\Delta t \cdot \mathcal{F}_{1}^{1:N_t}$}
\draw[blue] (6.5,0.25) -- (6.5,1.25) -- (8.5,1.25) -- (8.5,0.25) -- (6.5,0.25);

\put(205,57){$-\Delta t \cdot \mathcal{F}_{1}^{0:N_{t}-1}$}
\draw[red] (6.5,1.5) -- (6.5,2.5) -- (8.5,2.5) -- (8.5,1.5) -- (6.5,1.5);

\draw [black, ->  ] (8.5,2) -- (9,2.0);

\draw [black, ->  ] (8.5,3.5) -- (9,2.1);

\draw [black, ->  ] (8.5,5) -- (9.1,2.3);

\draw [black, ->  ] (8.5,1) -- (9,1.9);

\draw [black, -  ] (7.5,-0.5) -- (7.5,0.0);

\draw [black, -  ] (7.5,0.0) -- (5.25,0);

\draw [black, -  ] (5.25,0.0) -- (5.25,1.4);

\draw [black, -  ] (5.25,1.4) -- (5.6,1.4);

\put(175,67){\color{red}$FWD$}
\draw [red, -  ] (5.6,1.4) -- (5.6,2);
\draw [red, -> ] (5.6,2.0) -- (6.5,2);

\put(175,30){\color{blue}$BWD$}
\draw [blue, -  ] (5.6,1.4) -- (5.6,0.8);
\draw [blue, -> ] (5.6,0.8) -- (6.5,0.8);

\draw[] (9.36,2) circle (12pt);
\put(288,60){$\sum$}

\put(308,68){$\mathcal{K}_{g}^{0:N_{t}-1}$}
\draw [black, ->  ] (9.75,2) -- (11.25,2);

\end{tikzpicture}
\caption{Example DC-RNN for determining the $g$-Equation \eqref{eq:GGG} based on Forward Euler (Red) and Backward Euler (Blue) schemes. The inputs are $\rho(t_{n})$ and $g(t_{n})$ for $n=0,1,2,\cdots, N_t$. The dictionary contains order $\mathcal{O}(1)$ and $\mathcal{O}(\veps)$ operators. These operators are generated by the RNNs corresponding to orders $\veps^{-m}$ $m=0,1$ using $\mathcal{A}_{1}$, $\mathcal{A}_{2}$. The output $\mathcal{K}_{g}^{n}$ ($n=0,1,\cdots,N_t-1$) is to be minimized with respect to a chosen norm.}
\end{figure}

\textbf{Forward Euler scheme.}
The forward Euler scheme only involves two time steps and, hence, $q=1$. We can specify $\mathcal{K}_{u}^{n}(D^{n,1};\bd{\theta})$ to relate the data pair
\begin{equation}
D^{n,1}=\{u(x,t_{n}),u(x,t_{n+1})\}
\end{equation}
using a forward finite difference approximation for $\partial_{t}u=\mathcal{F}(u(x,t);\bm{\theta})$, the right hand side of which is a symbolic RNN as an equation ansatz. This gives us the forward Euler fitting scheme:
\begin{equation}\label{eq:fwrkappa}
\mathcal{K}_{u}^{n}(D^{n,1};\bd{\theta})=u(x,t_{n+1})-u(x,t_{n})-\Delta t\cdot \mathcal{F}(u(x,t_{n});\bm{\theta}).
\end{equation}
Minimizing the loss in Equation \eqref{eq:unreg} will determine a PDE governing the training data with time accuracy $\Delta t$. We display in Figure \ref{NNfigFwrd} a DC-RNN for determining the equation satisfied by $g(v,x,t)$ based on the Forward Euler scheme.

\textbf{Backward Euler scheme.}
The Backward Euler scheme for $\partial_{t}u=\mathcal{F}(u(x,t);\bm{\theta})$ relates the data pair
\begin{equation}
D^{n,1}=\{u(x,t_{n}),u(x,t_{n+1})\}
\end{equation}
using the backward Euler fitting scheme:
\begin{equation}\label{eq:bwrkappa}
\mathcal{K}_{u}^{n}(D^{n,1};\bd{\theta})=u(x,t_{n+1})-u(x,t_{n})-\Delta t\cdot \mathcal{F}(u(x,t_{n+1});\bm{\theta}).
\end{equation}
Minimizing the loss in Equation \eqref{eq:unreg} will determine a PDE governing the training data with accuracy $\Delta t$. We display in Figure \ref{NNfigFwrd} a DC-RNN for determining the equation satisfied by $g(v,x,t)$ based on the Backward-Euler scheme.

We only focused on the scalar equation in \eqref{eq:PTexpansionMain} to illustrate the loss function for the above schemes. The construction for the two-component vector equation in \eqref{eq:vectorizedansatz} is similar. The loss function is the sum of the loss function for each component
\begin{equation}\label{eq:unreg2}
L(\bd{\theta})=\dfrac{1}{N_{t}-q}\sum_{n=1}^{N_t-q}||\mathcal{K}_{g}^{n}(D^{n,q};\bd{\theta})||_*+||\mathcal{K}_{\rho}^{n}(D^{n,q};\bd{\theta})||_*,
\end{equation}
where $\mathcal{K}_{g}^{n}(D^{n,q};\bd{\theta})$ and $\mathcal{K}_{\rho}^{n}(D^{n,q};\bd{\theta})$ relate the data in
\begin{equation*}
D^{n,q}:=\{g(x,t_{n}),g(x,t_{n+1}),\cdots,g(x,t_{n+q}),\rho(x,t_{n}),\rho(x,t_{n+1}),\cdots,\rho(x,t_{n+q})\},\hspace{0.5cm}n=1,\cdots N_{t}-q.
\end{equation*}
 
\textbf{First-Order IMEX Scheme.}
In this paper, we are interested in Equations \eqref{eq:rho} and \eqref{eq:GGG} and, hence, will use schemes specialized for them. 

A first-order IMEX scheme for solving Equations \eqref{eq:rho} and \eqref{eq:GGG} are given by
\begin{equation}\label{eq:IMEX1G}
\begin{split}
g_{i+1/2}^{n+1}&=g_{i+1/2}^{n}+\Delta t\left\{\dfrac{1}{\veps}\left(I-\langle\,\, \rangle\right)\left(v^{+}\dfrac{g_{i+1/2}^{n}-g_{i-1/2}^{n}}{\Delta x}+v^{-}\dfrac{g_{i+3/2}^{n}-g_{i+1/2}^{n}}{\Delta x}\right)\right.\\
&\left.-\dfrac{\sigma_{i+1/2}^{S}}{\veps^{2}}g_{i+1/2}^{n+1}-\dfrac{1}{\veps^{2}}v\dfrac{\rho_{i+1}^{n}-\rho_{i}^{n}}{\Delta x}-\sigma_{i+1/2}^{A}g_{i+1/2}^{n}\right\},
\end{split}
\end{equation}
\begin{equation}\label{eq:IMEX1rho}
\rho_{i}^{n+1}=\rho_{i}^{n}+\Delta t\left\{\left\langle v\dfrac{g_{i+1/2}^{n+1}-g_{i-1/2}^{n+1}}{\Delta x}\right\rangle-\sigma_{i}^{A}\rho_{i}^{n}+G_{i}\right\},
\end{equation}
where $v^{+}=\dfrac{v+|v|}{2}$ and $v^{-}=\dfrac{v-|v|}{2}$. From this, we see that Equation \eqref{eq:IMEX1G} gives a relationship among $D^{n,1}$ that can be generalized to the ansatz:

\begin{equation}\label{eq:imex1fit}
\begin{split}
\mathcal{K}_{g}^{n}(D^{n,1};\bd{\theta})=&(1+\Delta t\dfrac{\sigma^{S}(x)}{\veps^{2}})g(v,x,t_{n+1})-(1-\Delta t\sigma^{A}(x))g(v,x,t_{n})\\
-&\Delta t \cdot \mathcal{F}_{1}(g(v,x,t_{n}),\rho(x,t_{n});\bd{\theta}),
\end{split}
\end{equation}
where $\mathcal{F}_{1}$ is the operator ansatz introduced in \eqref{eq:vectorizedansatz}.

Equation \eqref{eq:IMEX1rho} gives a relationship between data in $D^{n,1}$ via:
\begin{equation}\label{eq:imex1fit3}
\begin{split}
\mathcal{K}_{\rho}^{n}(D^{n,1};\bd{\theta})&=\rho(x,t_{n+1})-(1-\Delta t \sigma_{A}(x))\rho(x,t_{n})-\Delta t G(x)\\
&-\Delta t \cdot \mathcal{F}_{2}(g(v,x,t_{n+1}),\rho(x,t_{n});\bd{\theta}).
\end{split}
\end{equation}

$\mathcal{F}_{1}$ and $\mathcal{F}_{2}$ will be learned by minimizing the loss function \eqref{eq:unreg2}. In fact, one does not need to assume that the functions $\sigma^{S}(x)$, $\sigma^{A}(x)$, and $G(x)$ are known. One can learn these functions during the training process by replacing them with neural networks. For the special case when $\sigma^{S}(x)$ and $\sigma^{A}(x)$ are constants, we can replace them with trainable parameters $w_{S}$ and $w_{A}$, respectively.

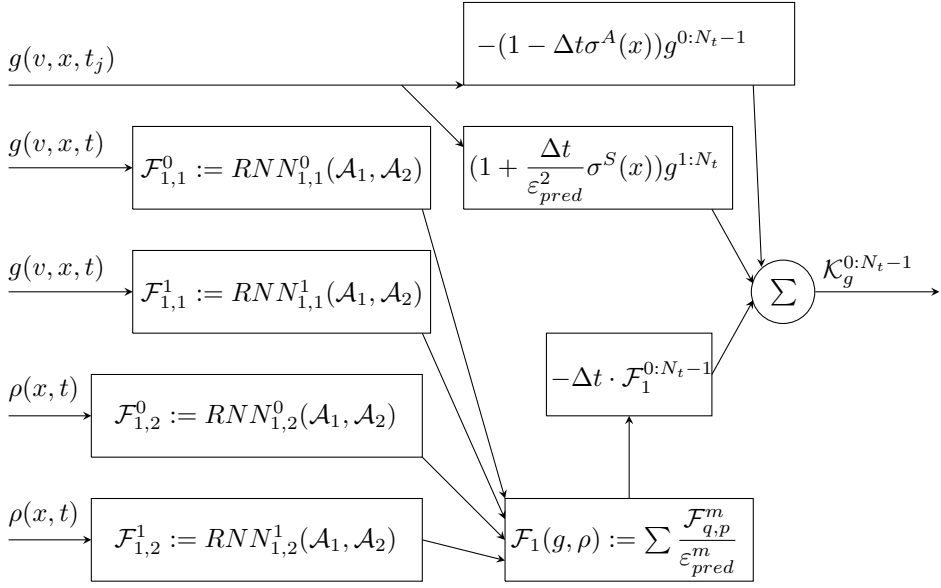
\begin{figure}\label{NNfigIMEXfirst}
\begin{tikzpicture}[x=1.1cm, y=1.1cm, >=stealth]

\draw [black, ->          ] (0,4.5) -- (5.5,4.5);
\put(0,148){$g(v,x,t_{j})$}

\draw [black, ->          ] (4.75,4.5) -- (5.5,3.75);

\draw [black, ->          ] (0,3.5) -- (1.5,3.5);
\put(0,116){$g(v,x,t)$}
\put(50,106){$\mathcal{F}_{1,1}^{0}:=RNN_{1,1}^{0}(\mathcal{A}_{1}, \mathcal{A}_{2})$}
\draw[black] (1.5,3) -- (1.5,4) -- (5.1,4) -- (5.1,3) -- (1.5,3);

\draw [black, ->          ] (0,2) -- (1.5,2);
\put(0,69){$g(v,x,t)$}
\put(50,59){$\mathcal{F}_{1,1}^{1}:=RNN_{1,1}^{1}(\mathcal{A}_{1}, \mathcal{A}_{2})$}
\draw[black] (1.5,1.5) -- (1.5,2.5) -- (5.1,2.5) -- (5.1,1.5) -- (1.5,1.5);

\draw [black, ->          ] (0,0.5) -- (1,0.5);
\put(0,23){$\rho(x,t)$}
\put(40,13){$\mathcal{F}_{1,2}^{0}:=RNN_{1,2}^{0}(\mathcal{A}_{1}, \mathcal{A}_{2})$}
\draw[black] (1,0) -- (1,1) -- (5,1) -- (5,0) -- (1,0);

\draw [black, ->          ] (0,-1) -- (1,-1);
\put(0,-24){$\rho(x,t)$}
\put(40,-34){$\mathcal{F}_{1,2}^{1}:=RNN_{1,2}^{1}(\mathcal{A}_{1}, \mathcal{A}_{2})$}
\draw[black] (1,-1.5) -- (1,-0.5) -- (5,-0.5) -- (5,-1.5) -- (1,-1.5);

\put(190,-34){$\mathcal{F}_{1}(g,\rho):=\sum\dfrac{\mathcal{F}_{q,p}^{m}}{\veps_{pred}^{m}}$}
\draw (6,-1.5) -- (6,-0.5) -- (9,-0.5) -- (9,-1.5) -- (6,-1.5);

\draw [black, ->  ] (5,3) -- (6,-0.5);

\draw [black, ->  ] (5,1.5) -- (6,-0.75);

\draw [black, ->  ] (5,0) -- (6,-1);

\draw [black, ->  ] (5,-1.0) -- (6,-1.25);

\put(177,153){$-(1-\Delta t\sigma^{A}(x))g^{0:N_t-1}$}
\draw (5.5,4.5) -- (5.5,5.5) -- (9.5,5.5) -- (9.5,4.5) -- (5.5,4.5);

\put(174,107){$(1+\dfrac{\Delta t}{\veps_{pred}^{2}} \sigma^{S}(x))g^{1:N_t}$}
\draw (5.5,3) -- (5.5,4) -- (8.75,4) -- (8.75,3) -- (5.5,3);

\put(205,27){$-\Delta t \cdot \mathcal{F}_{1}^{0:N_t-1}$}
\draw (6.5,0.5) -- (6.5,1.5) -- (8.5,1.5) -- (8.5,0.5) -- (6.5,0.5);

\draw [black, ->  ] (9,4.5) -- (9.1,2.3);

\draw [black, ->  ] (8.5,3) -- (9,2.1);

\draw [black, ->  ] (8.5,1) -- (9,1.9);

\draw [black, ->  ] (7.5,-0.5) -- (7.5,0.5);

\draw[] (9.36,2) circle (12pt);
\put(288,60){$\sum$}

\put(308,68){$\mathcal{K}_{g}^{0:N_{t}-1}$}
\draw [black, ->  ] (9.75,2) -- (11.25,2);

\end{tikzpicture}
\caption{Example DC-RNN for determining the $g$-Equation \eqref{eq:GGG} based on the First order IMEX scheme. The inputs are $\rho(t_{n})$ and $g(t_{n})$ for $n=0,1,2,\cdots,N_t$. The dictionary contains order $\mathcal{O}(1)$ and $\mathcal{O}(\veps)$ operators. These operators are generated by the RNNs corresponding to orders $\veps^{-m}$ $m=0,1$ using $\mathcal{A}_{1}$, $\mathcal{A}_{2}$. The output $\mathcal{K}_{g}^{n}$ ($n=0,1,\cdots,N_t-1$) is to be minimized with respect to a chosen norm.}
\end{figure}

We display in Figure \ref{NNfigIMEXfirst} a DC-RNN for determining the equation satisfied by $g(v,x,t)$ based on the First order IMEX scheme. One can go higher order with high-order IMEX schemes. These will either introduce more intermediate stages (if using IMEX Runge-Kutta schemes) or relate more data points to each other by increasing $q$ (if using IMEX multistep schemes). We leave details concerning high-order IMEX schemes in the Appendix section.

\subsection{Constructing $\langle g\rangle=0$ into our DC-RNN}

To improve physically accurate predictions, we recommend doing prior analysis of the data. In case the residuals for $g$ suggests that
\begin{equation}
\langle g\rangle=\int_{-1}^{1}g\,dv =0,
\end{equation}
the network can be further improved with a few modifications. First observe the following about the projection operator:

{\lemma \label{lemma_reg_0}
Define $\mathcal{A}_{p}:=\int_{-1}^{1}g(x,v,t)\,dv$ to be the projection operator, if $\mathcal{A}_{p}(g)=0$, then $\mathcal{A}_{i_{1}}\circ\cdots\circ\mathcal{A}_{i_m}\circ\mathcal{A}_{p}(g)=0$ for any sequence of linear operators $\mathcal{A}_{i_{k}}$.}

{\proof This is straightforward since if $\mathcal{A}_{p}$ is applied first to $g$, then the following linear operators will act on $0$.}

We restrict the lemma to linear operators since nonlinear operators such as $\mathcal{A}:=\exp(\cdot)$ may lead to a nonzero result. For example, in this case $\mathcal{A}\circ\mathcal{A}_{p}(g)=1$.

To apply lemma \ref{lemma_reg_0} to the DC-RNN in section \ref{dictionarybuilding}, we need to omit the identity and projection operators in $\mathcal{C}_{2}$ in equation \eqref{eq:ope1}. Equivalently, one can set the weights for these terms to be $0$. 

\remark The condition $\langle g\rangle=0$ is sometimes difficult for a machine learning algorithm to learn via regularization. For this condition to hold,
\begin{equation}\label{eq:showconstraint}
\langle \mathcal{F}_{1}(g,\rho)\rangle=0,
\end{equation}
as can be justified for the Forward-Euler case by performing the following calculation:
\begin{equation}
\begin{split}
\langle g(v,x,t_{n+1}) \rangle&\approx\langle g(v,x,t_{n})+\Delta t \cdot \mathcal{F}_{1}\rangle\,\,\,(\text{using Forward Euler})\\
&=\Delta t\cdot \langle \mathcal{F}_{1}\rangle\,\,\,\,\,\,\,\,\,\,\,\,\,\,\,\,\,\,\,\,\,\,\,\,\,\,\,\,\,\,\,(\text{by linearity and }\,\,\langle g\rangle=0).
\end{split}
\end{equation}
Thus, $\langle g\rangle=0$ implies \eqref{eq:showconstraint}. The justification of \eqref{eq:showconstraint} for other Runge-Kutta schemes is similar.

The loss enforcing Equation \eqref{eq:showconstraint} is:
\begin{equation}\label{eq:const_with_lagrange_mult}
L(\bd{\theta})=\dfrac{1}{N_t-q}\sum_{n=1}^{N_{t}-q}||\mathcal{K}_{g}^{n}(D;\bd{\theta})||+||\Delta t\cdot \langle \mathcal{F}_{1}(g,\rho;\bd{\theta}) \rangle||.
\end{equation}

We recommend using this only for black-box networks and not symbolic neural networks such as the DC-RNN.

\textbf{Regularization via Sparsity.} If the size of the dictionary is too large, one is more likely to over-fit data to an incorrect PDE. To help omit terms that do not appear in the PDE, we impose sparsity in weights and biases in the DC-RNN. Denote the set of all trainable parameters excluding $w_{eps}$, the parameter which trains $\veps_{pred}$, by $\bd{\theta}$. The regularization term
\begin{equation}\label{eq:reg_of_sparce}
R(\bd{\theta})=\gamma\cdot||\bd{\theta}||_{\ell^{1}}\,\,\text{with}\,\,\gamma\in\mathbb{R}^{+}
\end{equation}
is one natural choice to produce a PDE with the fewest possible terms. $\gamma$ is typically chosen to be a small number. In our numerical examples, $\gamma$ is chosen to be of order $10^{-4}$.

\textbf{Regularization via the Continuity of Weights and Biases.} We note that if $\sigma^{S}(x)$, $\sigma^{A}(x)$, $G(x)$, weights, or biases are not constant in $x$, we need to introduce neural networks to parametrize them to capture the dynamics of these functions of $x$. In this paper, any trainable function of $x$ will be fitted similar to a piecewise n-th degree polynomial function of x. Similar to spline interpolation, we will promote continuity of the $n-1$ derivatives. The polynomial basis, $p_{i}(x)$, $\{i=0,\cdots,N_p-1\}$ will be small with $N_p$ typically one order of magnitude smaller than $N_x$. Each $p_{i}(x)$ will be defined on the i-th partition of $x\in [a,b]$,
\begin{equation}
[a_{i},b_{i}]=[a+i\Delta x,a+(i+1)\Delta x],\,\,\Delta x=\frac{(b-a)}{N_{p}},\,\, i=0,1,\cdots,N_{p}-1
\end{equation}
Additional parameters $\{a_{i,j}\}$ are introduced because of the basis set
\begin{equation}
p_{i}(x):=a_{i,0}+a_{i,1}x+\cdots+a_{i,n}x^{n},
\end{equation}
enlarging our parameter set $\bd{\theta}$.
To promote the continuity in $x$ and its gradients at $a_{i}$,  $\{i=1,\cdots,N_{p}-1\}$, we apply the regularization term $R(\bd{\theta})$ to the loss function:
\begin{equation}\label{eq:cont}
R(\bd{\theta})=\gamma\sum_{i=1}^{deg(p_{i}(x))-1}\sum_{j=0}^{N_{p}-2}\cdot|\nabla_{x}^{i}p_{j}(a_{j+1})-\nabla_{x}^{i}p_{j+1}(a_{j+1})|\,\,\text{with}\,\,\gamma\in\mathbb{R}^{+}.
\end{equation}
Physically, Equation \eqref{eq:cont} is used to lessen the jump discontinuity in the learned functions of $x$. 

\section{Optimization}\label{sec:opt}
In this section we discuss how we update parameters $(\veps_{pred},\bd{\theta})$ to reach a minimum for our loss function. Since we do not assume that $\veps$ is known a priori, we will have to train this parameter. At the end of this section, we argue why our algorithm is expected to be superior to existing algorithms.
Our reasoning suggests that terms of order $\mathcal{O}(\veps^{n})$ should be updated with a learning rate proportional to ${\veps^{n}}$. We will verify this claim through several numerical experiments in the next section.

\subsection{Training of \texorpdfstring{$\veps_{pred}$}{Lg}}

A major goal for our algorithm is to determine approximately the magnitude of the scale $\veps$ involved in the multiscale dynamics. To fulfill the condition $0< \veps\leq 1$, the $\veps_{pred}$ in our algorithm is set to
\begin{equation}\label{eq:nonparallel}
\veps_{pred}=\dfrac{1}{2}(\text{tanh}(w_{\veps})+1),
\end{equation}
where $w_{\veps}$ is a trainable parameter. However, if one can parallelize, it makes more sense to restrict $\veps_{pred}$ over several intervals spanning $(0,1]$. For instance, let $s(i)= 0.1^i$ and
\begin{equation}\label{eq:parallel}
\veps_{pred}=\dfrac{s(i)-s(i+1)}{2}(\text{tanh}(w_{\veps}^{i})+\text{min}_{i}),
\end{equation}
where $(0,1]=[s(1),s(2)]\cup[s(2),s(3)]\cup[s(3),s(4)]\cup\cdots$. After training over each interval, one can choose the PDE corresponding to the lowest loss. Thus, with Equation \eqref{eq:parallel}, one has better control over where local minimums of the loss function occur.

\subsection{Training of parameters}

The parameters for our loss function \eqref{eq:unreg}, can be trained using our suggested algorithm: Adam method \cite{Kingma2014AdamAM}. This algorithm is great at training a relatively large number of parameters efficiently. Other gradient descent methods are possible including stochastic gradient descent. We will discuss an implementation of stochastic gradient descent below using the loss equation \eqref{eq:unreg} with respect to the $L_{1}$ and $L_{2}$ norms. While not necessary, we will simplify the calculations by using the forward Euler approximation and assuming $R(\theta)$ involves only sparse regularity. We can rewrite equation \eqref{eq:unreg} as:
\begin{equation}\label{eq:LplusR1}
\wh{L}(w_{\veps},\bd{\theta},\bd{x})=L(w_{\veps},\bd{\theta},\bd{x})+R_{1}
\end{equation}
where the regularization is
\begin{equation}
\begin{split}
R_{1}&:=\gamma_{1}||\bd{\theta}||_{\ell^{1}}\\
&=\gamma_{1}\sum_{i}|\theta_{i}|
\end{split}
\end{equation}
and the loss function is given by
\begin{equation}
L(w_{\veps},\bd{\theta},\bd{x})=\dfrac{1}{N_{t}}\sum_{j=1}^{N_{t}}||u(\bd{x},t_{j+1})-u(\bd{x},t_{j})+\int_{t_{j}}^{t_{j}+\Delta t}\sum_{n}\dfrac{1}{\veps(w_{\veps})^{n}} \mathcal{F}^{n}(u(\bd{x},s),\bd{\theta}_{n})\,ds||_{\ast}.\\
\end{equation}
Using the Forward Euler approximation,
\begin{equation}
\begin{split}
L(w_{\veps},\bd{\theta},\bd{x})&=L_{Fwrd}(w_{\veps},\bd{\theta},\bd{x})\\
&=\dfrac{1}{N_{t}}\sum_{j=1}^{N_{t}}||\mathcal{K}_{u}^{j}(u(\bd{x},t_{j}),\bd{\theta})||_{\ast}\\
&=\dfrac{1}{N_{t}}\sum_{j=1}^{N_{t}}||u(\bd{x},t_{j+1})-u(\bd{x},t_{j})+\Delta t\sum_{n}\dfrac{1}{\veps(w_{\veps})^{n}} \mathcal{F}^{n}(u(\bd{x},t_{j}),\bd{\theta}_{n})||_{\ast}
\end{split}
\end{equation}
where $(w_{\veps},\bd{\theta}):=(w_{\veps},\bd{\theta}_{0},\bd{\theta}_{1},\cdots)$. The gradient of equation \eqref{eq:LplusR1} with respect to $L_{1}$ and $L_{2}$ is given by
\begin{equation}\label{eq:gradients}
\begin{split}
\nabla_{(w_{\veps},\bd{\theta})}\wh{L}&=\dfrac{1}{N_{t}}\sum_{j=1}^{N_{t}}\text{sign}(L)(\partial_{w_{\veps}}L,\Delta t\nabla_{\bd{\theta}_{0}}\mathcal{F}^{0},\Delta t\dfrac{1}{\veps(w_{\veps})}\nabla_{\bd{\theta}_{1}}\mathcal{F}^{1},\cdots)+\gamma_{1}\text{sign}(\bd{\theta})\hspace{1cm}\text{using}\,\,L^{1}\\
\nabla_{(w_{\veps},\bd{\theta})}\wh{L}&=\dfrac{1}{N_{t}}\sum_{j=1}^{N_{t}}2L(\partial_{w_{\veps}}L,\Delta t\nabla_{\bd{\theta}_{0}}\mathcal{F}^{0},\Delta t\dfrac{1}{\veps(w_{\veps})}\nabla_{\bd{\theta}_{1}}\mathcal{F}^{1},\cdots)+\gamma_{1}\text{sign}(\bd{\theta})\hspace{1cm}\text{using}\,\,L^{2}
\end{split}
\end{equation}

If one records many data (large $N_{t}$), the summation in equation \eqref{eq:gradients} can be slow to compute. Thus, one can reduce computational resources by using our suggested \textbf{stochastic gradient descent}: The sum is taken over a random smaller subset of $\{1,2,\cdots,N_{t}\}$ of size $N_{s}<N_{t}$ and we replace $N_{t}$ with $N_{s}$ in equation \eqref{eq:gradients}.

\subsection{Discussion}
We will now discuss the effect of asymptotic expansion and our sparse regularization method. Because it is difficult to obtain clean algebraic expressions using ``all of" $\wh{L}$, we will consider a quadratic Taylor series truncation of $\mathcal{K}:=\mathcal{K}_{u}^{j}(u(\bd{x},t_{j}),\bd{\theta})$ near the point that minimizes $\wh{L}$ which we denote by  $(w_{\veps}^{\ast},\bd{\theta}_{0}^{\ast},\bd{\theta}_{1}^{\ast},\cdots)$. We will again assume the forward Euler approximation for $\mathcal{K}$. The point $(w_{\veps}^{(0)},\bd{\theta}_{0}^{(0)},\bd{\theta}_{1}^{(0)},\cdots)$ will denote the initial values for the training parameters and $(w_{\veps}^{(k)},\bd{\theta}_{0}^{(k)},\bd{\theta}_{1}^{(k)},\cdots)$ will denote the $k$-th step taken by the gradient descent process. The gradient at the $k$-th step is given by,
\begin{equation}
\begin{split}
\bd{z}:=&\nabla_{(w_{\veps},\bd{\theta})}\mathcal{K}_{u}^{j}(w_{\veps}^{(k)},\bd{\theta}^{(k)})\\
=&(\partial_{w_{\veps}}\mathcal{K}_{u}^{j},\Delta t\nabla_{\bd{\theta}_{0}}\mathcal{F}^{0},\Delta t\dfrac{1}{\veps(w_{\veps})}\nabla_{\bd{\theta}_{1}}\mathcal{F}^{1},\cdots,\Delta t\dfrac{1}{{\veps(w_{\veps})}^{M}}\nabla_{\bd{\theta}_{M}}\mathcal{F}^{M})|_{(w_{\veps}^{(k)},\bd{\theta}_{0}^{(k)},\bd{\theta}_{1}^{(k)},\cdots,\bd{\theta}_{M}^{(k)})}\\
=&(\partial_{w_{\veps}}\mathcal{K}_{u}^{j},\wt{\bd{z}})\\
=&(\partial_{w_{\veps}}\mathcal{K}_{u}^{j},\wt{\bd{z}}_{0},\wt{\bd{z}}_{1},\cdots,\wt{\bd{z}}_{M})
\end{split}
\end{equation}
where we defined
\begin{equation}
\begin{split}
\wt{\bd{z}}&:=(\wt{\bd{z}}_{0},\wt{\bd{z}}_{1},\cdots,\wt{\bd{z}}_{M})\\
&:=(\Delta t\nabla_{\bd{\theta}_{0}}\mathcal{F}^{0},\Delta t\dfrac{1}{\veps(w_{\veps})}\nabla_{\bd{\theta}_{1}}\mathcal{F}^{1},\cdots,\Delta t\dfrac{1}{{\veps(w_{\veps})}^{M}}\nabla_{\bd{\theta}_{M}}\mathcal{F}^{M})
\end{split}
\end{equation}
in order to simplify the notation. The hessian is given by

\begin{equation}
\begin{split}
H:&=\nabla^{2}\mathcal{K}_{u}^{j}(w_{\veps}^{(k)},\bd{\theta}^{(k)})\\
&=\Delta t\left[\begin{matrix}
\dfrac{1}{\Delta t}\partial_{w_{\veps}}^{2}\mathcal{K}_{u}^{j} & [\bd{0}]_{1\times d} & -\dfrac{1}{\veps^{2}}\partial_{w_{\veps}}\veps\nabla_{\bd{\theta}_{1}}\mathcal{F}^{1} & -\dfrac{2}{\veps^{3}}\partial_{w_{\veps}}\veps\nabla_{\bd{\theta}_{2}}\mathcal{F}^{2} & \cdots\\
\nabla_{\bd{\theta}_{0}}(\partial_{w_{\veps}}\mathcal{K}_{u}^{j})^{T} & [\nabla_{\bd{\theta}_{0}}^{2}\mathcal{F}^{0}]_{d\times d} & [0]_{d\times d} & [0]_{d\times d} &\\
\nabla_{\bd{\theta}_{1}}(\partial_{w_{\veps}}\mathcal{K}_{u}^{j})^{T} & [0]_{d\times d} & [\dfrac{1}{\veps}\nabla_{\bd{\theta}_{1}}^{2}\mathcal{F}^{1}]_{d\times d} & [0]_{d\times d} &\\
\nabla_{\bd{\theta}_{2}}(\partial_{w_{\veps}}\mathcal{K}_{u}^{j})^{T} &  [0]_{d\times d} &  [0]_{d\times d} &  [\dfrac{1}{\veps^{2}}\nabla_{\bd{\theta}_{2}}^{2}\mathcal{F}^{2}]_{d\times d} & \\
\vdots & & & & \ddots
\end{matrix}\right]
\end{split}
\end{equation}
where $d$ is the dimension of each $\bd{\theta}_{n}$. We also denote the $(d\cdot M-1)\times (d\cdot M-1)$ submatrix of $H$ by
\begin{equation}
\wt{H}=\Delta t\left[\begin{matrix}
[\nabla_{\bd{\theta}_{0}}^{2}\mathcal{F}^{0}]_{d\times d} & [0]_{d\times d} & [0]_{d\times d} &\\
[0]_{d\times d} & [\dfrac{1}{\veps}\nabla_{\bd{\theta}_{1}}^{2}\mathcal{F}^{1}]_{d\times d} & [0]_{d\times d} &\\
[0]_{d\times d} &  [0]_{d\times d} &  [\dfrac{1}{\veps^{2}}\nabla_{\bd{\theta}_{2}}^{2}\mathcal{F}^{2}]_{d\times d} & \\
\vdots & & & \ddots
\end{matrix}\right]
\end{equation}

\textbf{The effect of asymptotic expansion on learning rate.}
We now consider the effect of updating the parameters $\bd{\theta}$ via gradient descent. To further simplify algebraic expressions, we will assume that $\veps(w_{\veps}^{\ast})=\veps^{\ast}=\veps$ is the constant optimal value. According to the gradient descent method:
\begin{equation}\label{eq:update}
\begin{split}
\bd{\theta}^{(k+1)}\leftarrow& \bd{\theta}^{(k)}-\bd{\alpha}\otimes\wt{\bd{z}}\\
:&=(\bd{\theta}^{(k)})-(\alpha_{0}\wt{\bd{z}}_{0},\alpha_{1}\wt{\bd{z}}_{1},\cdots,\alpha_{M}\wt{\bd{z}}_{M})
\end{split}
\end{equation}
where the parameter $\bd{\alpha}:=(\alpha_{0},\alpha_{1},\alpha_{2},\cdots)$ will denote the learning rate which updates step $k\rightarrow k+1$. $\alpha_{n}$ will denote the learning rate for the terms of order $\mathcal{O}(\veps^{-n})$. Our goal is to understand the optimal behaviour of the learning rates $\alpha_{n}$. Substituting equation \eqref{eq:update} into our quadratic truncation of $\mathcal{K}_{n}^{j}$ yields:

\begin{equation}\label{eq:thetaupdate1}
\begin{split}
\mathcal{K}_{n}^{j}(w_{\veps}^{\ast},\bd{\theta}^{(k+1)})&=\mathcal{K}_{n}^{j}(w_{\veps}^{\ast},\bd{\theta}^{(k)})-(\bd{\alpha}\otimes\wt{\bd{z}})^{T}\wt{\bd{z}}+(\bd{\alpha}\otimes\wt{\bd{z}})^{T}\wt{H}(\bd{\alpha}\otimes\wt{\bd{z}})
\end{split}
\end{equation}

when $(\bd{\alpha}\otimes\wt{\bd{z}})^{T}H(\bd{\alpha}\otimes\wt{\bd{z}})$ is positive we can solve for the optimal values for $\bd{\alpha}$:
\begin{equation}\label{eq:thetaupdate2}
\bd{\alpha}=\dfrac{\wt{\bd{z}}^{T}\wt{\bd{z}}}{\wt{\bd{z}}^{T}\wt{H}\wt{\bd{z}}}\hspace{0.5cm} \iff \hspace{0.5cm} \alpha_{n}=\dfrac{{\veps^{\ast}}^{n}\wt{\bd{z}}_{n}^{T}\wt{\bd{z}}_{n}}{\Delta t\wt{\bd{z}}_{n}^{T}[\nabla_{\bd{\theta}_{n}}^{2}\mathcal{F}^{n}]_{d\times d}\wt{\bd{z}}_{n}}
\end{equation}
The meaning of the above calculations is summarize below:

\textbf{Observation.} If $\mathcal{F}^{n}$ is well approximated by a quadratic function with $\veps(w_{\veps})=\eps^{\ast}$, then $\bd{\theta}_{n}^{(k+1)}\leftarrow \bd{\theta}_{n}^{(k)}$ should be updated (according to equations \eqref{eq:thetaupdate1} and \eqref{eq:thetaupdate2}) in the direction of  $\nabla_{\theta_{n}}\mathcal{F}^{n}$. The optimal learning rate is proportional to $\dfrac{1}{\veps^{n}}$. The eigenvalues and vectors of $[\nabla_{\bd{\theta}_{n}}^{2}\mathcal{F}^{n}]_{d\times d}$ determine the stability of the learning process. In the worst case scenario, $\wt{\bd{z}}_{n}$ is in the direction corresponding to the largest eigenvector of $[\nabla_{\bd{\theta}_{n}}^{2}\mathcal{F}^{n}]_{d\times d}$.

\begin{example}

Because the observation above made use of several simplifying assumptions, we provide some numerical evidence to support this claim. What we observe through repeated numerical tests is that the multiscale fitting methods tend to converge to the correct model using less training time and iterations. Evidence of this is shown in figure \ref{Loss_comparison}. We also note that convergence to a lower loss value does not necessarily translate to a better prediction.


\begin{figure}[ht]
\begin{center}
\centerline{\includegraphics[width=3.5in]{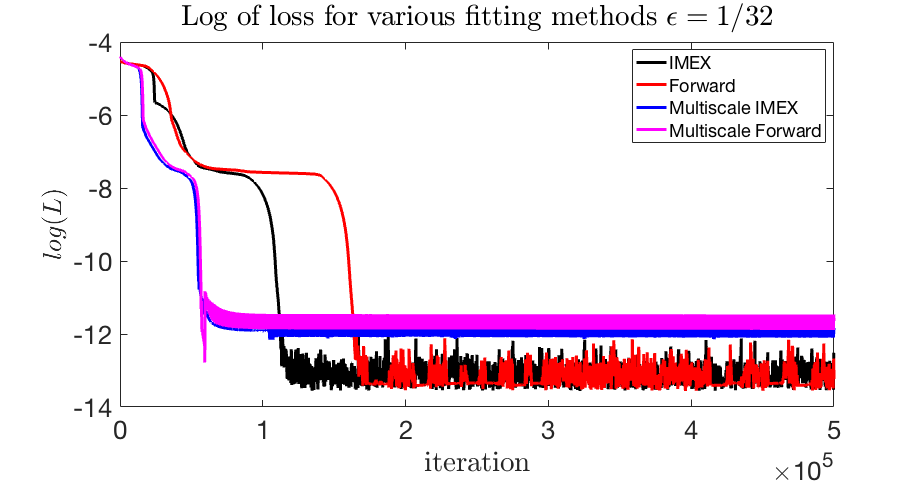}\includegraphics[width=3.5in]{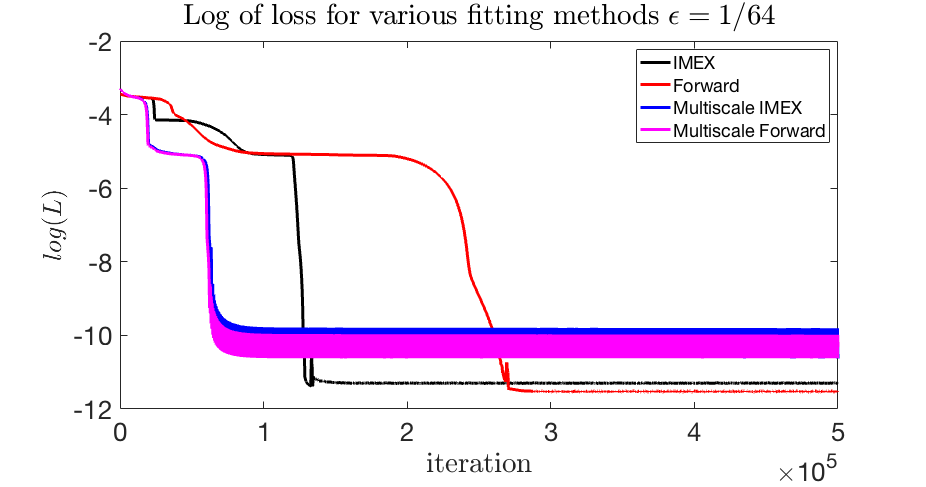}}
\centerline{\includegraphics[width=3.5in]{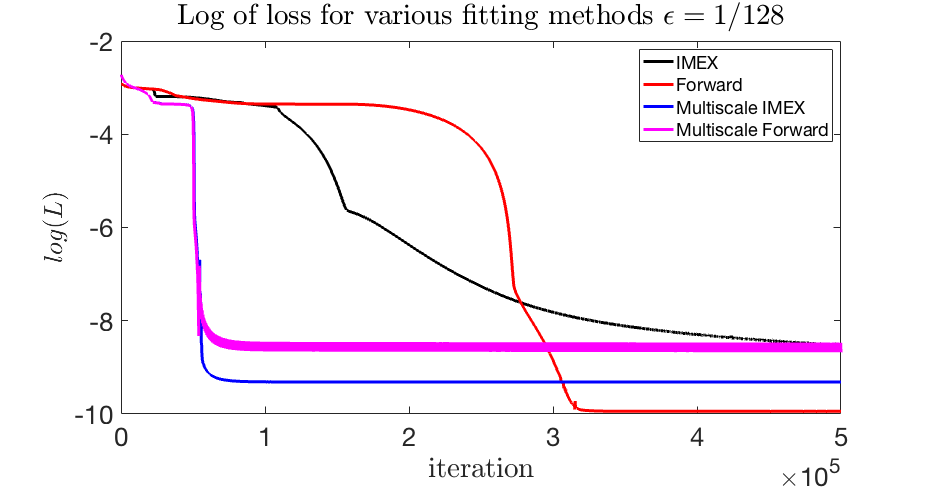}\includegraphics[width=3.5in]{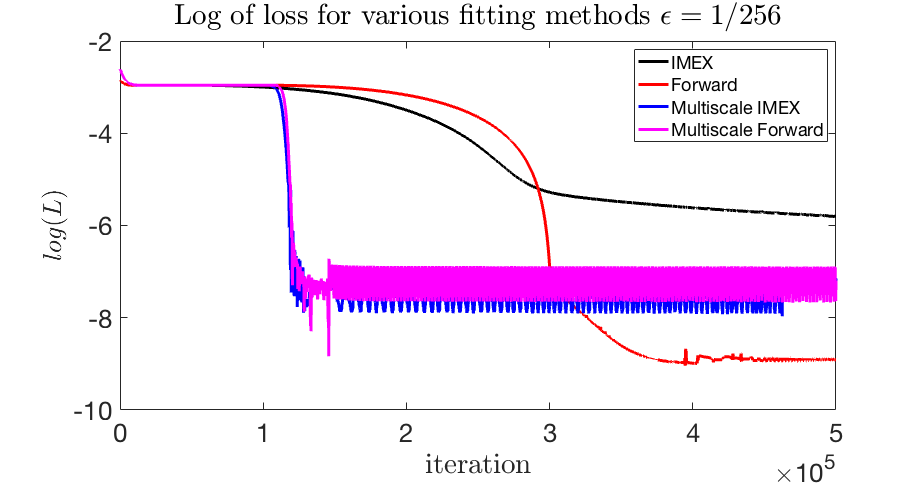}}
\caption{The behaviour of the loss function for example \ref{first_example} is displayed above. We see quick convergence to the a small loss when using the multiscale fitting methods, though a very small loss might mean overfitting. As we can see in other tables in Section \ref{sec:results}, the estimation accuracy of IMEX is better than forward Euler and the accuracy with multiscale ansatz is better than that without multiscale ansatz.}
\label{Loss_comparison}
\end{center}
\vskip -0.2in
\end{figure}

\end{example}

\remark We would like to remark that the efficiency of the proposed DC-RNN is demonstrated experimentally. As we made some simplifying assumptions on the terms $\mathcal{F}^{n}(\bd{\theta}_{n})$ generated by our RNNs, we acknowledge that theoretical analysis remains vastly open, though several seminal works have been available \cite{AllenZhu2019CanSL,NEURIPS2019_0ee8b85a,keller2020discovery,yu2020onsagernet,li2020curse}.

\remark other algorithms such as \cite{Liu2020HierarchicalDL,BinD17,DBLP:journals/corr/abs-1711-10561,DBLP:journals/corr/abs-1711-10566,RAISSI2018125,articleCKBSK,wang2018deep}, do not have an adaptive
$\veps_{pred}$. Numerically, having an adaptive $\veps_{pred}$ mimics adaptive gradient descent methods. The behavior of the loss functions vs. iteration, has typically converged with respect to the number of iterations when compared to using the Adam algorithm with no multiscale expansion.

\textbf{The effect of sparse regularization.} The conventional Lasso method uses sparse regularization where one fits the data to an ansatz with an $L_{1}$ penalty as in equation \eqref{eq:LplusR1}. Typically the coefficients for the basis set of functions in the hypothesis space is set to be a sparse vector. Using our notation, this means
\begin{equation}\label{eq:l1_conventional}
R_{1}:=\gamma_{1}\sum_{m\geq 1}\sum_{\pi\in\mathcal{D}}|a_{\pi(1),\cdots,\pi(m)}(\bd{\theta})|\hspace{1cm}(\text{for some}\,\,\gamma_{1}\in\mathbb{R}^{+})
\end{equation}
is added to the loss function \eqref{eq:intro}. However, in the case that some of these coefficients are large $\mathcal{O}(\veps^{-n})$ for $n>0$ and $0<\veps\ll 1$, equation \eqref{eq:l1_conventional} will create problems for the learning process. The issue is that there will be conflicting goals: keeping a particular set of coefficients $a_{\pi(1),\cdots,\pi(m)}(\bd{\theta})$ large in magnitude while at the same time minimizing $R_{1}$ as much as possible. Even if one sets $\gamma_{1}$ to be a very small value, one will still run into the trouble of setting appropriate learning rates as mentioned earlier. This conflict is clearly solved by our algorithm using the regularization
\begin{equation}\label{eq:l1_nonconventional}
R_{1}:=\gamma_{1}\sum_{m\geq 1}\sum_{\pi\in\mathcal{D}}|a_{\pi(1),\cdots,\pi(m)}(\bd{\theta}_{0})|+|a_{\pi(1),\cdots,\pi(m)}(\bd{\theta}_{1})|+\cdots+|a_{\pi(1),\cdots,\pi(m)}(\bd{\theta}_{M})|
\end{equation}
and setting the coefficients of the basis terms to
\begin{equation}
a_{\pi(1),\cdots,\pi(m)}(\bd{\theta}_{0})+\dfrac{a_{\pi(1),\cdots,\pi(m)}(\bd{\theta}_{1})}{\veps}+\cdots+\dfrac{a_{\pi(1),\cdots,\pi(m)}(\bd{\theta}_{M})}{\veps^{M}}.
\end{equation}
With this design we are able to have both sparsity and large $\mathcal{O}(\veps^{-n})$ coefficients. The only drawback is that we had to introduce more parameters for our design.

\section{Numerical Examples}\label{sec:results}

\begin{table}[t]
\label{Tab:forwardEulermethods}
\vskip 0.15in
\begin{center}
\begin{small}
\begin{sc}
\begin{tabular}{lccr}
\toprule
$\veps$ & Multiscale & learned $g$-equation using Foward Euler Scheme & Error Type-I/ II \\
\midrule
1/ 16  & No & $\partial_{t}g=-( 16 ^2+{\color{red} 0.836 \cdot 10^{ -1 }}) g -( 16 +{\color{red} 2.052 \cdot 10 ^{ -1 }}) v\cdot \partial_{x}g$  & \\
& & $+ ( 16 +{\color{red} 2.066 \cdot 10^{ -1 }})\langle v\partial_{x}g\rangle -( 16^{2} -{\color{red} 1.145 \cdot 10^{ -1 }})v\cdot\partial_{x}\rho+\cdots$ & { 0.11 \%}\,/\,0.74\%\\
1/ 32  & No & $\partial_{t}g=-( 32 ^2-{\color{red}1.855 }) g -( 32 +{\color{red} 3.220 \cdot 10 ^{ -1 }}) v\cdot \partial_{x}g$  & \\
& & $+ ( 32 +{\color{red} 3.272 \cdot 10^{ -1 }})\langle v\partial_{x}g\rangle -( 32^{2} -{\color{red} 2.061 })v\cdot\partial_{x}\rho+\cdots$ & { 0.21\%}\,/\,0.60\%\\
1/ 64  & No & $\partial_{t}g=-( 64 ^2-{\color{red} 3.302 \cdot 10^{ 1 }}) g -( 64 +{\color{red} 1.389 \cdot 10 ^{ -1 }}) v\cdot \partial_{x}g$  & \\
& & $+ ( 64 +{\color{red} 1.321 \cdot 10^{ -1 }})\langle v\partial_{x}g\rangle -( 64^{2} -{\color{red} 3.317 \cdot 10^{ 1 }})v\cdot\partial_{x}\rho+\cdots$ & { 0.79 \%}\,/\,0.33\%\\
1/ 128  & No & $\partial_{t}g=-( 128 ^2-{\color{red} 0.524 \cdot 10^{ 3 }}) g -( 128 -{\color{red} 3.512}) v\cdot \partial_{x}g$  & \\
& & $+ ( 128 -{\color{red} 3.541})\langle v\partial_{x}g\rangle -( 16384 -{\color{red} 0.524 \cdot 10^{ 3 }})v\cdot\partial_{x}\rho+\cdots$ & { 3.19 \%}\,/\,2.98\%\\
1/ 256  & No & $\partial_{t}g=-( 256 ^2-{\color{red} 0.788 \cdot 10^{ 4 }}) g -( 256 -{\color{red} 3.059 \cdot 10 ^{ 1 }}) v\cdot \partial_{x}g$  & \\
& & $+ ( 256 -{\color{red} 3.066 \cdot 10^{ 1 }})\langle v\partial_{x}g\rangle -( 256^2 -{\color{red} 0.788 \cdot 10^{ 4 }})v\cdot\partial_{x}\rho+\cdots$ & { 12.03 \%}\,/\,11.99\%\\
\bottomrule
\end{tabular}
\end{sc}
\end{small}
\end{center}\caption{Learned $g$-equation using the DC-RNN algorithm based on Forward-Euler scheme without multiscale ansatz.}
\vskip -0.1in
\end{table}

\begin{table}
\label{Tab:IMEXmethods}
\vskip 0.15in
\begin{center}
\begin{small}
\begin{sc}
\begin{tabular}{lccr}
\toprule
$\veps$ & Multiscale & learned $g$-equation using IMEX1 Scheme & Error Type-I/ II \\
\midrule
1/ 16  & No & $\partial_{t}g=-( 16 ^2+{\color{red} 3.344 \cdot 10^{ -1 }}) g -( 16 +{\color{red} 2.210 \cdot 10 ^{ -1 }}) v\cdot \partial_{x}g$  & \\
& & $+ ( 16 +{\color{red} 2.227 \cdot 10^{ -1 }})\langle v\partial_{x}g\rangle -( 16^{2} +{\color{red} 1.542 \cdot 10^{ -1 }})v\cdot\partial_{x}\rho+\cdots$ & { 0.17 \%}\,/\,0.74\%\\
1/ 32  & No & $\partial_{t}g=-( 32 ^2-{\color{red} 1.319}) g -( 32 +{\color{red} 3.478 \cdot 10 ^{ -1 }}) v\cdot \partial_{x}g$  & \\
& & $+ ( 32 +{\color{red} 3.539 \cdot 10^{ -1 }})\langle v\partial_{x}g\rangle -( 32^{2} -{\color{red} 1.538})v\cdot\partial_{x}\rho+\cdots$ & { 0.16 \%}\,/\,0.62\%\\
1/ 64  & No & $\partial_{t}g=-( 64 ^2+{\color{red} 2.291 \cdot 10^{ 1 }}) g -( 64 +{\color{red} 1.000}) v\cdot \partial_{x}g$  & \\
& & $+ ( 64 +{\color{red} 1.004})\langle v\partial_{x}g\rangle -( 64^{2} +{\color{red} 2.224 \cdot 10^{ 1 }})v\cdot\partial_{x}\rho+\cdots$ & { 0.56 \%}\,/\,1.06\%\\
1/ 128  & No & $\partial_{t}g=-( 128 ^2+{\color{red} 4.745 \cdot 10^{ 2 }}) g -( 128 +{\color{red} 4.348 }) v\cdot \partial_{x}g$  & \\
& & $+ ( 128 +{\color{red} 4.335 })\langle v\partial_{x}g\rangle -( 128^{2} +{\color{red} 4.755 \cdot 10^{ 2 }})v\cdot\partial_{x}\rho+\cdots$ & { 2.90 \%}\,/\,3.15\%\\
1/ 256  & No & $\partial_{t}g=-( 256 ^2+{\color{red} 2.915 \cdot 10^{ 3 }}) g -( 256 +{\color{red} 1.1598 \cdot 10 ^{ 1 }}) v\cdot \partial_{x}g$  & \\
& & $+ ( 256 +{\color{red} 1.156 \cdot 10^{ 1 }})\langle v\partial_{x}g\rangle -( 256^2 -{\color{red} 2.931 \cdot 10^{ 3 }})v\cdot\partial_{x}\rho+\cdots$ & { 4.46 \%}\,/\,4.49\%\\
\bottomrule
\end{tabular}
\end{sc}
\end{small}
\end{center}\caption{Learned $g$-equation using the DC-RNN algorithm based on IMEX1 scheme without multiscale ansatz.}
\vskip -0.1in
\end{table}

In this section, we test our DC-RNN using the PDE example in \eqref{eq:rho} and \eqref{eq:GGG} with various values of $\veps$. In the numeral results presented in this section, to make it intuitive to clarify the estimation error, the predicted coefficients are of the form
\begin{equation}
\text{predicted} = (\text{exact} + \text{{\color{red}difference}}),
\end{equation}
where in red we highlight the difference from the predicted to the exact value. The smaller the magnitude of the difference, the better the prediction. We also include the percentage error of Type-I defined via the relative error of the PDE coefficients in the $\ell^1$-norm:
\begin{equation}
\text{Error Type-I}:=\dfrac{\sum_{i}|\text{exact coefficients}_{i} -\text{predicted coefficients}_{i}|}{\sum_{i}|\text{exact coefficients}_{i}|}\times 100\%.
\end{equation}
To achieve fair comparisons, we will also use the percentage error of Type-II defined as the average relative errors of the prediction of each nonzero PDE coefficient:
\begin{equation}
    \text{Error Type-II}:=\dfrac{1}{\#|\text{nonzero terms}|}\sum_{i}\dfrac{|\text{exact coefficient}_{i}-\text{predicted coefficient}_{i}|}{|\text{exact coefficient}|_{i}}\times 100\%.
\end{equation}
Obviously when $\text{exact coefficient}_{i}=0$, the sum is undefined, thus we will only sum over the non-zero coefficients of the PDE.

The right hand side of the learned PDE will contain many terms, for the sake of readability, we display only the terms involved in either Equation \eqref{eq:rho} or \eqref{eq:GGG} in our numerical results. Other terms are typically minute in magnitude due to our sparsity regularization.

\textbf{Data Gathering.} The data that we produce in our examples are computed with IMEX-ARS(2,2,2) scheme using small mesh size $\Delta x=\frac{1}{1000}$ and $\Delta t=\frac{1}{2}\Delta x^{2}$. Thus, the data can be assumed to be nearly an exact solution to Equations \eqref{eq:rho} and \eqref{eq:GGG}. The velocity discretization we used in our examples is the standard 16-point Gauss quadrature in $[-1,1]$.


The training data is prepared by taking a subset of the exact data to reduce the memory cost. We define the training set number of grid points by $\wt{N}_{x}$ and $\wt{N}_{t}$ and note that for all examples $\wt{N}_{v}=N_{v}=16$. To obtain a subset of the data points, a coarser grid is chosen: $\wt{\Delta x}\geq \Delta x$ and $\wt{\Delta t}\geq \Delta t$ with $\wt{\Delta x}$ and $\wt{\Delta t}$ satisfying $\Delta t=\mathcal{O}(\Delta x^2)$. Code will be made available at \url{https://github.com/Ricard0000} or the authors' personal homepages.

\begin{example}\label{first_example}
\textbf{Forward Euler vs. IMEX1 and Multiscale vs. Non-multiscale.}
In this experiment, we compare the fitting using the forward Euler and the first order IMEX scheme. We ran our algorithm using the RNN \eqref{eq:eq210} with one layer. We will also compare the results with those using a multiscale ansatz ($M=2$ in (\ref{eq:PTexpansionMain})). The data was produced with $\sigma^{S}(x)=1$, $\sigma^{A}(x)=0$. We choose $\wt{N}_{x}=1000$ and $\wt{N}_{t}=56$.
We only attempt to learn the dynamics of the $g$-Equation \eqref{eq:GGG} to shorten the length of this paper. The results are shown in Tables \ref{Tab:forwardEulermethods}, \ref{Tab:IMEXmethods}, \ref{Tab:forwardEulermethodsMS}, and \ref{Tab:IMEXmethodsMS}, where we can see that the IMEX1 can produce a relatively better result compared with forward Euler, especially when $\varepsilon$ is small, no matter the multiscale ansatz is used or not. The fitting using the multiscale ansatz, on the other hand, does not show a clear advantage or disadvantage in this example, because, in low-order methods, the discretization error dominates the benefits of the multiscale ansatz. This is the reason why we will study high-order methods to demonstrate the advantage of the multiscale ansatz in the next example.



\begin{table}[t]
\label{Tab:forwardEulermethodsMS}
\vskip 0.15in
\begin{center}
\begin{small}
\begin{sc}
\begin{tabular}{lccr}
\toprule
$\veps$ & Multiscale & learned $g$-equation using Forward Euler Scheme & Error Type-I/ II \\
\midrule
1/ 16  & Yes & $\partial_{t}g=-( 16 ^2-{\color{red} 3.616 \cdot 10^{ -2 }}) g -( 16 +{\color{red} 2.161 \cdot 10 ^{ -1 }}) v\cdot \partial_{x}g$  & \\
& & $+ ( 16 +{\color{red} 2.116 \cdot 10^{ -1 }})\langle v\partial_{x}g\rangle -( 16^2 -{\color{red} 0.866 \cdot 10^{ -1 }})v\cdot\partial_{x}\rho+\cdots$ & { 0.10 \%}\,/\,0.08\%\\
1/ 32  & Yes & $\partial_{t}g=-( 32 ^2-{\color{red} 2.113}) g -( 32 +{\color{red} 3.833 \cdot 10 ^{ -1 }}) v\cdot \partial_{x}g$  & \\
& & $+ ( 32 +{\color{red} 3.623 \cdot 10^{ -1 }})\langle v\partial_{x}g\rangle -( 32^{2} -{\color{red} 1.882})v\cdot\partial_{x}\rho+\cdots$ & { 0.22 \%}\,/\,0.68\%\\
1/ 64  & Yes & $\partial_{t}g=-( 64 ^2-{\color{red} 3.325 \cdot 10^{ 1 }}) g -( 64 +{\color{red} 1.774 \cdot 10 ^{ -1 }}) v\cdot \partial_{x}g$  & \\
& & $+ ( 64 +{\color{red} 1.887 \cdot 10^{ -1 }})\langle v\partial_{x}g\rangle -( 64^{2} -{\color{red} 3.414 \cdot 10^{ 1 }})v\cdot\partial_{x}\rho+\cdots$ & { 0.81 \%}\,/\,0.55\%\\
1/ 128  & Yes & $\partial_{t}g=-( 128 ^2-{\color{red} 0.522 \cdot 10^{ 3 }}) g -( 128 -{\color{red} 3.225}) v\cdot \partial_{x}g$  & \\
& & $+ ( 128 -{\color{red} 3.259 })\langle v\partial_{x}g\rangle -( 128^2 -{\color{red} 0.526 \cdot 10^{ 3 }})v\cdot\partial_{x}\rho+\cdots$ & { 3.19 \%}\,/\,2.87\\
1/ 256  & Yes & $\partial_{t}g=-( 256 ^2+{\color{red} 0.787 \cdot 10^{ 4 }}) g -( 256 -{\color{red} 2.884 \cdot 10 ^{ 1 }}) v\cdot \partial_{x}g$  & \\
& & $+ ( 256 -{\color{red} 3.067 \cdot 10^{ 1 }})\langle v\partial_{x}g\rangle -( 256^2 +{\color{red} 0.787 \cdot 10^{ 4 }})v\cdot\partial_{x}\rho+\cdots$ & { 12.01 \%}\,/\,11.82\%\\
\bottomrule
\end{tabular}
\end{sc}
\end{small}
\end{center}\caption{Learned $g$-equation using the DC-RNN algorithm based on Forward-Euler scheme with multiscale ansatz.}
\vskip -0.1in
\end{table}

\begin{table}[t]
\label{Tab:IMEXmethodsMS}
\vskip 0.15in
\begin{center}
\begin{small}
\begin{sc}
\begin{tabular}{lccr}
\toprule
$\veps$ & Multiscale & learned $g$-equation using IMEX1 Scheme & Error Type-I/ II \\
\midrule
1/ 16  & Yes & $\partial_{t}g=-( 16 ^2+{\color{red} 0.617 \cdot 10^{ -2 }}) g -( 16 +{\color{red} 2.033 \cdot 10 ^{ -1 }}) v\cdot \partial_{x}g$  & \\
& & $+ ( 16 +{\color{red} 2.094 \cdot 10^{ -1 }})\langle v\partial_{x}g\rangle -( 16^{2} -{\color{red} 1.282 \cdot 10^{ -1 }})v\cdot\partial_{x}\rho+\cdots$ & { 0.10 \%}\,/\,0.64\%\\
1/ 32  & Yes & $\partial_{t}g=-( 32 ^2-{\color{red} 1.863 }) g -( 32 +{\color{red} 3.303 \cdot 10 ^{ -1 }}) v\cdot \partial_{x}g$  & \\
& & $+ ( 32 +{\color{red} 3.466 \cdot 10^{ -1 }})\langle v\partial_{x}g\rangle -( 32^{2} -{\color{red} 1.826})v\cdot\partial_{x}\rho+\cdots$ & { 0.20 \%}\,/\,0.62\%\\
1/ 64  & Yes & $\partial_{t}g=-( 64 ^2-{\color{red} 2.672 \cdot 10^{ 1 }}) g -( 64 +{\color{red} 2.679 \cdot 10 ^{ -1 }}) v\cdot \partial_{x}g$  & \\
& & $+ ( 64 +{\color{red} 2.887 \cdot 10^{ -1 }})\langle v\partial_{x}g\rangle -( 64^2 -{\color{red} 2.758 \cdot 10^{ 1 }})v\cdot\partial_{x}\rho+\cdots$ & { 0.65 \%}\,/\,0.55\%\\
1/ 128  & Yes & $\partial_{t}g=-( 128 ^2-{\color{red} 2.488 \cdot 10^{ 2 }}) g -( 128 -{\color{red} 1.139 }) v\cdot \partial_{x}g$  & \\
& & $+ ( 128 -{\color{red} 1.011})\langle v\partial_{x}g\rangle -( 128^2 -{\color{red} 2.523 \cdot 10^{ 2 }})v\cdot\partial_{x}\rho+\cdots$ & { 1.52 \%}\,/\,1.19\%\\
1/ 256  & Yes & $\partial_{t}g=-( 256 ^2-{\color{red} 0.586 \cdot 10^{ 4 }}) g -( 256 -{\color{red} 2.537 \cdot 10 ^{ 1 }}) v\cdot \partial_{x}g$  & \\
& & $+ ( 256 -{\color{red} 2.156 \cdot 10^{ 1 }})\langle v\partial_{x}g\rangle -( 256^2 -{\color{red} 0.586 \cdot 10^{ 4 }})v\cdot\partial_{x}\rho+\cdots$ & { 8.94 \%}\,/\,9.05\%\\
\bottomrule
\end{tabular}
\end{sc}
\end{small}
\end{center}\caption{Learned $g$-equation using the DC-RNN algorithm based on IMEX1 scheme with multiscale ansatz.}
\vskip -0.1in
\end{table}
\end{example}

\begin{example}
\textbf{Higher Order Methods: Multiscale vs. Non-multiscale.} In this experiment, we set $\sigma^{S}(x)=1$ and $\sigma^{A}(x)=0$ and use two layers in our RNN. We choose $\wt{N}_{x}=1000$ and $\wt{N}_{t}=56$ similarly to the previous example. Tests are done using the IMEX-BDF2 scheme and IMEX-ARS(2,2,2) scheme without and with the multiscale ansatz ($M=2$ in (\ref{eq:PTexpansionMain})). Results are recorded in Tables \ref{Tab:table1a}, \ref{Tab:table1b}, \ref{Tab:table1c}, and \ref{Tab:table1d}. The observation is as follows: High order IMEX methods, being either the BDF (multi-step) or ARS (multi-stage) can produce more accurate results than first order methods. However, even if using high-order methods, the prediction is getting worse for smaller $\varepsilon$. Then with the multiscale ansatz, good accuracy can be restored for small $\varepsilon$. All in all, high order methods combined with a multiscale ansatz produce stable and accurate results for a wide range of $\varepsilon$.


\end{example}

\begin{table}[t]
\label{Tab:table1a}
\vskip 0.15in
\begin{center}
\begin{small}
\begin{sc}
\begin{tabular}{lccr}
\toprule
$\varepsilon$  & Multiscale & learned $g$-equation using IMEX-BDF2 scheme & Error Type-I/ II \\
\midrule
1/ 16  & No & $\partial_{t}g=-( 16 ^2+{\color{red} 3.073 \cdot 10^{ -1 }}) g -( 16 -{\color{red} 0.531 \cdot 10 ^{ -1 }}) v\cdot \partial_{x}g$  & \\
& & $+ ( 16 -{\color{red} 0.715 \cdot 10^{ -3 }})\langle v\partial_{x}g\rangle -( 16^2 +{\color{red} 1.354 \cdot 10^{ -2 }})v\cdot\partial_{x}\rho+\cdots$ & { 0.06 \%}\,/\,0.11\%\\
1/ 32  & No & $\partial_{t}g=-( 32 ^2+{\color{red} 1.736 \cdot 10^{ -1 }}) g -( 32 +{\color{red} 0.615 \cdot 10 ^{ -1 }}) v\cdot \partial_{x}g$  & \\
& & $+ ( 32 +{\color{red} 1.243 \cdot 10^{ -2 }})\langle v\partial_{x}g\rangle -( 32^2 +{\color{red} 2.441 \cdot 10^{ -2 }})v\cdot\partial_{x}\rho+\cdots$ & { 0.01 \%}\,/\,0.06\%\\
1/ 64  & No & $\partial_{t}g=-( 64 ^2+{\color{red} 0.784}) g -( 64 -{\color{red} 0.777 \cdot 10 ^{ -2 }}) v\cdot \partial_{x}g$  & \\
& & $+ ( 64 -{\color{red} 1.195 \cdot 10^{ -2 }})\langle v\partial_{x}g\rangle -( 64^2 -{\color{red} 2.922 \cdot 10^{ -1 }})v\cdot\partial_{x}\rho+\cdots$ & { 0.01 \%}\,/\,0.01\%\\
1/ 128  & No & $\partial_{t}g=-( 128 ^2-{\color{red} 2.857 \cdot 10^{ 1 }}) g -( 128 -{\color{red} 0.953 \cdot 10 ^{ -1 }}) v\cdot \partial_{x}g$  & \\
& & $+ ( 128 -{\color{red} 2.176 \cdot 10^{ -1 }})\langle v\partial_{x}g\rangle -( 128^2 -{\color{red} 2.392 \cdot 10^{ 1 }})v\cdot\partial_{x}\rho+\cdots$ & { 0.15 \%}\,/\,0.14\%\\

1/ 256  & No & $\partial_{t}g=-( 256 ^2-{\color{red} 2.336 \cdot 10^{ 3 }}) g -( 256 -{\color{red} 2.436 \cdot 10 ^{ 1 }}) v\cdot \partial_{x}g$  & \\
& & $+ ( 256 -{\color{red} 0.901 \cdot 10^{ 1 }})\langle v\partial_{x}g\rangle -( 256^{2} -{\color{red} 2.334 \cdot 10^{ 3 }})v\cdot\partial_{x}\rho+\cdots$ & { 3.57 \%}\,/\,5.04\%\\

\bottomrule
\end{tabular}
\end{sc}
\end{small}
\end{center}\caption{Learned $g$-equation using the DC-RNN algorithm based on IMEX-BDF2 scheme without multiscale ansatz.}
\vskip -0.1in
\end{table}

\begin{table}[t]
\label{Tab:table1b}
\vskip 0.15in
\begin{center}
\begin{small}
\begin{sc}
\begin{tabular}{lccr}
\toprule
$\veps$ & Multiscale & learned $g$-equation using IMEX-ARS(2,2,2) scheme& Error Type-I/ II \\
\midrule
1/ 16  & No & $\partial_{t}g=-( 16 ^2+{\color{red} 1.707 \cdot 10^{ -1 }}) g -( 16 +{\color{red} 2.886 \cdot 10 ^{ -2 }}) v\cdot \partial_{x}g$  & \\
& & $+ ( 16 +{\color{red} 3.391 \cdot 10^{ -3 }})\langle v\partial_{x}g\rangle -( 16^2 -{\color{red} 1.251 \cdot 10^{ -3 }})v\cdot\partial_{x}\rho+\cdots$ & { 0.03 \%}\,/\,0.06\%\\
1/ 32  & No & $\partial_{t}g=-( 32 ^2+{\color{red} 3.706 \cdot 10^{ -1 }}) g -( 32 -{\color{red} 0.822 \cdot 10 ^{ -2 }}) v\cdot \partial_{x}g$  & \\
& & $+ ( 32 -{\color{red} 4.989 \cdot 10^{ -3 }})\langle v\partial_{x}g\rangle -( 32^2 +{\color{red} 3.662 \cdot 10^{ -2 }})v\cdot\partial_{x}\rho+\cdots$ & { 0.01 \%}\,/\,0.07\%\\

1/ 64  & No & $\partial_{t}g=-( 64 ^2-{\color{red} 0.689 \cdot 10^{ 1 }}) g -( 64 +{\color{red} 0.967 \cdot 10 ^{ -1 }}) v\cdot \partial_{x}g$  & \\
& & $+ ( 64 +{\color{red} 4.945 \cdot 10^{ -2 }})\langle v\partial_{x}g\rangle -( 64^2 +{\color{red} 4.291 \cdot 10^{ -1 }})v\cdot\partial_{x}\rho+\cdots$ & { 0.09 \%}\,/\,0.10\%\\

1/ 128  & No & $\partial_{t}g=-( 128 ^2-{\color{red} 0.602 \cdot 10^{ 3 }}) g -( 128 -{\color{red} 2.158 \cdot 10 ^{ -1 }}) v\cdot \partial_{x}g$  & \\
& & $+ ( 128 -{\color{red} 1.587 \cdot 10^{ -1 }})\langle v\partial_{x}g\rangle -(128^2 -{\color{red} 2.508 })v\cdot\partial_{x}\rho+\cdots$ & { 1.83 \%}\,/\,0.10\%\\

1/ 256  & No & $\partial_{t}g=-( 256 ^2-{\color{red} 0.814 \cdot 10^{ 3 }}) g -( 256 -{\color{red} 3.151 \cdot 10 ^{ -1 }}) v\cdot \partial_{x}g$  & \\
& & $+ ( 256 -{\color{red} 0.601})\langle v\partial_{x}g\rangle -( 256^2 -{\color{red} 1.880 \cdot 10^{ 2 }})v\cdot\partial_{x}\rho+\cdots$ & { 0.76\%}\,/\,0.47\%\\

\bottomrule
\end{tabular}
\end{sc}
\end{small}
\end{center}\caption{Learned $g$-equation using the DC-RNN algorithm based on IMEX-ARS(2,2,2) scheme without multiscale ansatz.}
\vskip -0.1in
\end{table}

\begin{table}[t]
\label{Tab:table1c}
\vskip 0.15in
\begin{center}
\begin{small}
\begin{sc}
\begin{tabular}{lccr}
\toprule
$\veps$ & Multiscale &learned $g$-equation using IMEX-BDF2 scheme & Error Type-I/ II \\
\midrule

1/ 16  & YES & $\partial_{t}g=-( 16 ^2-{\color{red} -1.217 \cdot 10^{ -1 }}) g -( 16 -{\color{red} -4.661 \cdot 10 ^{ -2 }}) v\cdot \partial_{x}g$  & \\
& & $+ ( 16 -{\color{red} 2.677 \cdot 10^{ -3 }})\langle v\partial_{x}g\rangle -( 256 -{\color{red} 3.619 \cdot 10^{ -2 }})v\cdot\partial_{x}\rho+\cdots$ & { 0.04 \%}\,/\,0.09\%\\

1/ 32  & YES & $\partial_{t}g=-( 32 ^2-{\color{red} 2.731 \cdot 10^{ -1 }}) g -( 32 -{\color{red} -0.972 \cdot 10 ^{ -1 }}) v\cdot \partial_{x}g$  & \\
& & $+ ( 32 -{\color{red} -1.634 \cdot 10^{ -2 }})\langle v\partial_{x}g\rangle -( 1024 -{\color{red} 4.102 \cdot 10^{ -1 }})v\cdot\partial_{x}\rho+\cdots$ & { 0.04 \%}\,/\,0.10\%\\

1/ 64  & YES & $\partial_{t}g=-( 64 ^2+{\color{red} 1.593} ) g -( 64 +{\color{red} 0.576 \cdot 10 ^{ -1 }}) v\cdot \partial_{x}g$  & \\
& & $+ ( 64 +{\color{red} 3.793 \cdot 10^{ -2 }})\langle v\partial_{x}g\rangle -( 64^2 +{\color{red} 1.709 })v\cdot\partial_{x}\rho+\cdots$ & { 0.04\%}\,/\,0.06\%\\
1/ 128  & YES & $\partial_{t}g=-( 128 ^2+{\color{red} 4.817 \cdot 10^{ 1 }}) g -( 128 -{\color{red} 0.573 \cdot 10 ^{ -1 }}) v\cdot \partial_{x}g$  & \\
& & $+ ( 128 +{\color{red} 1.349 \cdot 10^{ -1 }})\langle v\partial_{x}g\rangle -( 128^2 +{\color{red} 0.506 \cdot 10^{ 2 }})v\cdot\partial_{x}\rho+\cdots$ & { 0.29 \%}\,/\,0.19\%\\
1/ 256  & YES & $\partial_{t}g=-( 256 ^2+{\color{red} 4.745 \cdot 10^{ 2 }}) g -( 256 +{\color{red} 1.118 \cdot 10 ^{ 1 }}) v\cdot \partial_{x}g$  & \\
& & $+ ( 256 +{\color{red} 1.851})\langle v\partial_{x}g\rangle -( 256^2 +{\color{red} 4.688 \cdot 10^{ 2 }})v\cdot\partial_{x}\rho+\cdots$ & { 0.72 \%}\,/\,1.63\%\\
\bottomrule
\end{tabular}
\end{sc}
\end{small}
\end{center}\caption{Learned $g$-equation using the DC-RNN algorithm based on IMEX-BDF2 scheme with multiscale ansatz. }
\vskip -0.1in
\end{table}

\begin{table}[t]
\label{Tab:table1d}
\vskip 0.15in
\begin{center}
\begin{small}
\begin{sc}
\begin{tabular}{lccr}
\toprule
$\veps$ & Multiscale & learned $g$-equation using IMEX-ARS(2,2,2) scheme & Error Type-I/ II \\
\midrule

1/ 16  & YES & $\partial_{t}g=-( 16 ^2+{\color{red} 3.612 \cdot 10^{ -1 }}) g -( 16 +{\color{red} 1.024 \cdot 10 ^{ -2 }}) v\cdot \partial_{x}g$  & \\
& & $+ ( 16 +{\color{red} 0.685 \cdot 10^{ -3 }})\langle v\partial_{x}g\rangle -( 256 -{\color{red} 0.768 \cdot 10^{ -1 }})v\cdot\partial_{x}\rho+\cdots$ & { 0.08\%}\,/\,0.06\%\\

1/ 32  & YES & $\partial_{t}g=-( 32 ^2-{\color{red} 1.251 }) g -( 32 -{\color{red}2.309 \cdot 10 ^{ -2 }}) v\cdot \partial_{x}g$  & \\
& & $+ ( 32 +{\color{red} 2.640 \cdot 10^{ -3 }})\langle v\partial_{x}g\rangle -( 32^2 -{\color{red} 4.427 \cdot 10^{ -1 }})v\cdot\partial_{x}\rho+\cdots$ & { 0.08 \%}\,/\,0.06\%\\
1/ 64  & YES & $\partial_{t}g=-( 64 ^2-{\color{red} 3.075 \cdot 10^{ -2 }}) g -( 64 -{\color{red} 0.866 \cdot 10 ^{ -1 }}) v\cdot \partial_{x}g$  & \\
& & $+ ( 64 -{\color{red} 1.748 \cdot 10^{ -1 }})\langle v\partial_{x}g\rangle -( 64^2 +{\color{red} 1.177 })v\cdot\partial_{x}\rho+\cdots$ & { 0.01 \%}\,/\,0.11\%\\
1/ 128  & YES & $\partial_{t}g=-( 128 ^2+{\color{red} 1.018 \cdot 10^{ 1 }}) g -( 128 +{\color{red} 1.068 }) v\cdot \partial_{x}g$  & \\
& & $+ ( 128 +{\color{red} 2.269 \cdot 10^{ -1 }})\langle v\partial_{x}g\rangle -( 128^2+{\color{red} 1.044 \cdot 10^{ 1 }})v\cdot\partial_{x}\rho+\cdots$ & { 0.06 \%}\,/\,0.28\%\\

1/ 256  & YES & $\partial_{t}g=-( 256 ^2+{\color{red} 3.372 \cdot 10^{ 2 }}) g -( 256 +{\color{red} 2.386 }) v\cdot \partial_{x}g$  & \\
& & $+ ( 256 +{\color{red} 4.680 \cdot 10^{ -1 }})\langle v\partial_{x}g\rangle -( 65536 +{\color{red} 3.296 \cdot 10^{ 2 }})v\cdot\partial_{x}\rho+\cdots$ & { 0.51 \%}\,/\,0.53\%\\

\bottomrule
\end{tabular}
\end{sc}
\end{small}
\end{center}\caption{Learned $g$-equation using the DC-RNN algorithm based on IMEX-ARS(2,2,2) scheme with multiscale ansatz. }
\vskip -0.1in
\end{table}

\begin{figure}[ht]
\vskip -0.1in
\begin{center}
\centerline{\includegraphics[width=2.75in]{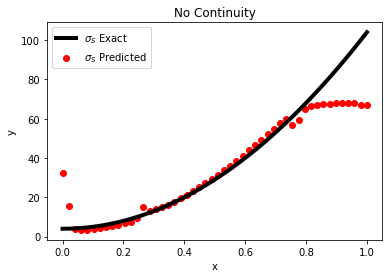}\includegraphics[width=2.75in]{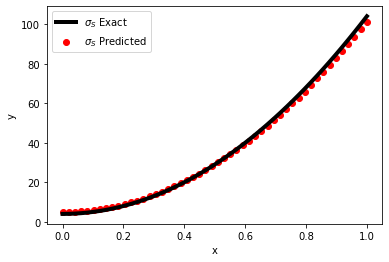}}
\caption{Left: Predicted $\sigma^{S}$ with no continuity constraint. Right: Predicted $\sigma^{S}$ with continuity constraint.}
\label{figsigma}
\end{center}
\vskip -0.2in
\end{figure}

\begin{example}
\textbf{Learning Space-Dependent Functions.}
We demonstrate that functions such as $\sigma^{S}(x)$, $\sigma^{A}(x)$, or $G(x)$ can be learned using space-dependent weights and biases. In this example, we choose
\begin{equation}
\sigma^{S}(x)=4+100x^{2},
\end{equation}
$\sigma^{A}(x)=0$, and $G(x)=0$. We set $\veps=1$ and use our DC-RNN based on IMEX fitting. The predicted PDE for the $g$-equation with no continuity regularization is:
\begin{equation}
\begin{split}
\partial_{t}g&=(1 +{\color{red}0.015})v\partial_{x}g -(1+{\color{red}0.017}) \langle v\partial_{x}g\rangle\\
 &+ (1+{\color{red}0.008} )v\partial_{x}\rho + [5-{\color{red}1.595} , 100-{\color{red}31.960}] g +\cdots,
\end{split}
\end{equation}
where $[4-{\color{red}0.595} , 100-{\color{red}31.960}]$ is the minimum and maximum values of $\sigma^{S}(x)$. We display the predicted $\sigma^{S}$ on the left of Figure \ref{figsigma}. We also impose the continuity regularization in Equation \eqref{eq:cont}, our predicted PDE is now given by:
\begin{equation}
\begin{split}
\partial_{t}g&=(1-{\color{red}0.004}) v\partial_{x}g -(1-{\color{red}0.005}) \langle v\partial_{x}g\rangle\\
&+ (1+{\color{red}0.001}) v\partial_{x}\rho + [ 5-{\color{red}0.048} , 100+{\color{red}}1.340] g +\cdots,
\end{split}
\end{equation}
with predicted $\sigma^{S}$ plotted on the right of Figure \ref{figsigma}. We note that the jump discontinuities on the left of Figure \ref{figsigma} are due to over fitting of the data. As we can see from this example,  these jumps have been removed by utilizing the continuity regularization in \eqref{eq:cont}.
\end{example}

\begin{example}
\textbf{The diffusion limit.}
As mentioned in Section 2, the equation for $\rho$ is given by \eqref{eq:rho}. However, when $\varepsilon \rightarrow 0$, it is well approximated by the diffusion equation \eqref{eq:diffusionlimit}. Therefore, we expect a good learning algorithm should be able to predict this diffusion limit when $\varepsilon$ is small. In this example, we will verify that the proposed DC-RNN is able to realize this diffusion limit. We choose $\sigma^{S}=1/3$, $\sigma^{A}=0$, and $G=0$. This means that in the limit $\veps\rightarrow 0$, the diffusion coefficient should be $1$ as an ideal test example for simplicity. We summarize the numerical results in Table \ref{Tab:diffusionexperiment}, where we see that the prediction indeed approaches to a diffusion equation with the diffusion coefficient equal to $1$.
\end{example}



\begin{table}[t]
\label{Tab:diffusionexperiment}
\vskip 0.15in
\begin{center}
\begin{small}
\begin{sc}
\begin{tabular}{lc}
\toprule
$\varepsilon$ & Learned $\rho$-equation \\
\midrule
$1/16$  & $\partial_{t}\rho=(-0.000141 ) \partial_{xx}\rho -(0.993171)\langle v\partial_{x}g\rangle+\cdots$  \\
\\
$1/256$ & $\partial_{t}\rho=(-0.0436861) \partial_{xx}\rho -(1.042619)\langle v\partial_{x}g\rangle +\cdots$ \\
\\
$1/2048$ & $\partial_{t}\rho=(0.985353) \partial_{xx}\rho -(0.003806)\langle v\partial_{x}g\rangle +\cdots$ \\
\\
$1/4096$ & $\partial_{t}\rho=(0.985596) \partial_{xx}\rho -(0.010862)\langle v\partial_{x}g\rangle+\cdots$ \\
\\
\bottomrule
\end{tabular}
\end{sc}
\end{small}
\end{center}\caption{Learned $\rho$-equations for various values of $\veps$ towards the diffusion limit when $\varepsilon\rightarrow 0$.}
\vskip -0.1in
\end{table}

\begin{example}\label{compare1}
\textbf{Comparison with Conventional Methods. (Part 1)}



The Lasso method \cite{Tibshirani1996RegressionSA,bookonlasso}  is a regression method used for variable selection. The method determines features based on $L_2$-minimization subject to sparse basis coefficients. This method does not assume a multiscale ansatz as in \eqref{eq:PTexpansionMain} or any other regularity conditions.  To apply the Lasso method for PDE discovery, we formulate a linear equation similar to \cite{Brunton3932} for the Lasso method to solve. The data $g$ and $\rho$ is mapped onto a dictionary of operators applied to the data and recorded in a matrix $A$, thus the memory requirements for using the Lasso method are typically larger compared to our algorithm. We also compute $\bd{b}:=\partial_{t}g$ using second order finite differences and allow the Lasso method to pick the best linear combination of the columns of the matrix $A$ that most closely resembles $\bd{b}$.

We perform tests of our algorithm vs. Lasso using $\sigma^{S}(x)=1$, $\sigma^{A}(x)=0$, and $G(x)=0$. We record results for the $g$-equation in Table \ref{Tab:Lasso}. The matrix $A$ in the Lasso method we implemented assumed 18 columns corresponding to terms involved in the dynamics (e.g., $g$, $v\partial_{x}\rho$, $v\partial_{x} g$, $\langle v\partial_{x}g\rangle$) and 14 others terms not actually involved  (e.g., artificially built by compositions of advection and projection operators \eqref{eq:theadvection} and \eqref{eq:theprojection}). We ran the Lasso method several times using several values of the sparsity regularization parameter $\alpha$. However, we only present the results associated with the best $\alpha$. The Lasso method performed fairly well but, it typically predicted more undesirable features for the dynamics and thus had a greater error than our DC-RNN as we can see from Table \ref{Tab:Lasso}.
\end{example}

\begin{table}[t]
\label{Tab:Lasso}
\vskip 0.15in
\begin{center}
\begin{small}
\begin{sc}
\begin{tabular}{lcccr}
\toprule
$\veps$ & Method & learned $g$-equation & Error/Relative Error \\
\midrule

1/ 16  & LASSO & $\partial_{t}g=-( 16 ^2-{\color{red} 1.598}) g -( 16 -{\color{red} 4.210 \cdot 10 ^{ -3 }}) v\cdot \partial_{x}g$  & \\
& & $+ ( 16 -{\color{red} 3.451 \cdot 10^{ -3 }})\langle v\partial_{x}g\rangle -( 16^2 -{\color{red} 4.715 \cdot 10^{ -2 }})v\cdot\partial_{x}\rho+\cdots$ & { 0.30 \%}\,/\,{0.17 \%}\\

1/ 32  & LASSO & $\partial_{t}g=-( 32 ^2-{\color{red} 4.512}) g -( 32 -{\color{red} 1.285 \cdot 10 ^{ -2 }}) v\cdot \partial_{x}g$  & \\
& & $+ ( 32 +{\color{red} 1.015})\langle v\partial_{x}g\rangle -( 32^2 -{\color{red} 0.521})v\cdot\partial_{x}\rho+\cdots$ & {0.28 \%}\,/\,{ 0.92 \%}\\

1/ 64  & LASSO & $\partial_{t}g=-( 64 ^2-{\color{red} 2.425 \cdot 10^{ 1 }}) g -( 64 -{\color{red} 1.291 \cdot 10 ^{ -1 }}) v\cdot \partial_{x}g$  & \\
& & $+ ( 64 +{\color{red} 0.735 \cdot 10^{ 1 }})\langle v\partial_{x}g\rangle -( 64^2 -{\color{red} 0.912 \cdot 10^{ 1 }})v\cdot\partial_{x}\rho+\cdots$ & { 0.491 \%}\,/\,{ 3.12 \%}\\

1/ 128  & LASSO & $\partial_{t}g=-( 128 ^2-{\color{red} 3.160 \cdot 10^{ 2 }}) g -( 128 -{\color{red} 2.175}) v\cdot \partial_{x}g$  & \\
& & $+ ( 128 +{\color{red} 3.725 \cdot 10^{ 1 }})\langle v\partial_{x}g\rangle -( 128^2 -{\color{red} 2.556 \cdot 10^{ 2 }})v\cdot\partial_{x}\rho+\cdots$ & { 1.85 \%}\,/\,{ 8.57 \%}\\

1/ 256  & LASSO & $\partial_{t}g=-( 256 ^2-{\color{red} 0.508 \cdot 10^{ 5 }}) g -( 256 -{\color{red} 2.162 \cdot 10 ^{ 2 }}) v\cdot \partial_{x}g$  & \\
& & $+ ( 256 -{\color{red} 1.740 \cdot 10^{ 2 }})\langle v\partial_{x}g\rangle -( 256 ^2-{\color{red} 0.505 \cdot 10^{ 5 }})v\cdot\partial_{x}\rho+\cdots$ & { 77.36 \%}\,/\,{ 76.79 \%}\\

\bottomrule
\end{tabular}
\end{sc}
\end{small}
\end{center}\caption{Learned $g$-equation using the Lasso method \cite{Tibshirani1996RegressionSA,bookonlasso} }
\vskip -0.1in
\end{table}

\begin{example}\label{compare2}
\textbf{Comparison with Conventional Methods. (Part 2)}

Next, we try the STRidge method in \cite{RSB16}. Similar to the Lasso method, a matrix of the dictionary is formed. Unlike the Lasso method, the STRidge method makes more efficient use of memory requirements and also features a hard threshold, i.e., large coefficients are assumed to be likely candidates for the dynamics of the PDE. Again, we use 18 terms for our dictionary as in the previous example. After running the STRidge algorithm, the predicted weights for the involved terms $g$, $v\partial_{x}\rho$, $v\partial_{x} g$, $\langle v\partial_{x}g\rangle$ are recorded in Table \ref{Tab:STridge}. The STRidge algorithm is very sensitive to noise, therefore it identified terms that are not supposed to be involved in the dynamics. The weights of the erroneous terms were so large that overall, the algorithm had a large error. For the STRidge algorithm, the main source of error is likely in the hard threshold assumption. 
\end{example}

\begin{table}[t]
\label{Tab:STridge}
\vskip 0.15in
\begin{center}
\begin{small}
\begin{sc}
\begin{tabular}{lccr}
\toprule
$\veps$ & Method & learned $g$-equation & Error/Relative Error \\
\midrule
1/ 16  & STRidge & $\partial_{t}g=-( 16 ^2+{\color{red} 1.318 \cdot 10^{ -1 }}) g -( 16 -{\color{red} 0.512 \cdot 10 ^{ -3 }}) v\cdot \partial_{x}g$  & \\
& & $+ ( 16 +{\color{red} 0.799 \cdot 10^{ 1 }})\langle v\partial_{x}g\rangle -( 16 ^2-{\color{red} 2.953 \cdot 10^{ -3 }})v\cdot\partial_{x}\rho+\cdots$ & { >99 \%}\,/\,{12.50 \%}\\

1/ 32  & STRidge & $\partial_{t}g=-( 32 ^2+{\color{red} 1.395 \cdot 10^{ -1 }}) g -( 32 +{\color{red} 4.808 \cdot 10 ^{ -4 }}) v\cdot \partial_{x}g$  &  \\
& & $+ ( 32 -{\color{red} 1.600 \cdot 10^{ 1 }})\langle v\partial_{x}g\rangle -( 32 ^2+{\color{red} 1.450 \cdot 10^{ -2 }})v\cdot\partial_{x}\rho+\cdots$ & { >99 \%}\,/\,{12.50 \%}\\

1/ 64  & STRidge & $\partial_{t}g=-( 64 ^2+{\color{red} 3.426 \cdot 10^{ -1 }}) g -( 64 +{\color{red} 3.832 \cdot 10 ^{ -3 }}) v\cdot \partial_{x}g$  &  \\
& & $+ ( 64 -{\color{red} 3.200 \cdot 10^{ 1 }})\langle v\partial_{x}g\rangle -( 64 ^2+{\color{red} 2.308 \cdot 10^{ -1 }})v\cdot\partial_{x}\rho+\cdots$ & { >99 \%}\,/\,{ 12.50\%}\\

1/ 128  & STRidge & $\partial_{t}g=-( 128 ^2+{\color{red} 1.201\cdot 10^{ 1 }}) g -( 128 +{\color{red} 2.101\cdot 10 ^{ -1 }}) v\cdot \partial_{x}g$  &  \\
& & $+ ( 128 +{\color{red} 0.639 \cdot 10^{ 2 }})\langle v\partial_{x}g\rangle -( 128 ^2-{\color{red} 1.196 \cdot 10^{ 1 }})v\cdot\partial_{x}\rho+\cdots$ & { >99 \%}\,/\,{ 14.85\%}\\

1/ 256  & STRidge & $\partial_{t}g=-( 256 ^2+{\color{red} 0.761 \cdot 10^{ 3 }}) g -( 256 +{\color{red} 0.723 \cdot 10 ^{ 1 }}) v\cdot \partial_{x}g$  & \\
& & $+ ( 256 -{\color{red} 1.265 \cdot 10^{ 2 }})\langle v\partial_{x}g\rangle -( 256 ^2+{\color{red} 0.761 \cdot 10^{ 3 }})v\cdot\partial_{x}\rho+\cdots$ & { >99 \%}\,/\,{ 13.64 \%}\\
\bottomrule
\end{tabular}
\end{sc}
\end{small}
\end{center}\caption{Learned $g$-equation using STRidge method \cite{RSB16}.}
\vskip -0.1in
\end{table}

\begin{example}
\textbf{Comparison with Conventional Methods. (Part 3)}

Now we discuss the Physics-Informed-Neural-Network (PINN) in \cite{DBLP:journals/corr/abs-1711-10561,DBLP:journals/corr/abs-1711-10566}. In \cite{DBLP:journals/corr/abs-1711-10561,DBLP:journals/corr/abs-1711-10566}, the authors suggest forming feed forward neural nets mapping the domain $(v,x,t)$ to the values of $g$ and $\rho$. Denote these networks as $\mathcal{N}_g$ and $\mathcal{N}_\rho$, respectively.
$\mathcal{N}_g$ and $\mathcal{N}_\rho$ can be considered as functions in $(v,x,t)$ and, hence, we can apply differential operators to $\mathcal{N}_g$ and $\mathcal{N}_\rho$. The loss which trains $\mathcal{N}_g$, $\mathcal{N}_\rho$, and the equations they satisfy is given by:
\begin{equation}\label{eq:MLoss}
||F_{1}(\mathcal{N}_\rho, \mathcal{N}_g)||_{L^2}^2+||F_{2}(\mathcal{N}_\rho, \mathcal{N}_g)||_{L^2}^2+||g-\mathcal{N}_{g}||_{L^2}^2+||\rho-\mathcal{N}_{\rho}||_{L^2}^2,
\end{equation}
where
\begin{equation}\label{eq:nondep}
F_1(\rho,g):=\partial_{t}g-\left(\lambda_{1}v\partial_{x}g+\lambda_{2}\langle v\partial_{x}g\rangle+\lambda_{3}v\partial_{x}\rho+\lambda_{4}g+\lambda_{5}g\right),
\end{equation}
and
\begin{equation}
F_2(\rho,g):=\partial_{t}\rho-\left(\lambda_{6}\langle v\partial_x g\rangle+\lambda_{7}\sigma^{A}\rho+\lambda_{8}G\right).
\end{equation}
The $L^2$ norm in the loss function is discretized using the training samples of $g$ and $\rho$. After minimizing the loss function over the network parameters and $\lambda_{1},\cdots,\lambda_{8}$, we can identify $\mathcal{N}_g$ and $\mathcal{N}_\rho$ fitting $g$ and $\rho$, respectively, and $F_1$ and $F_2$ specifying the governing equation of $g$ and $\rho$. Though this method is very powerful in many applications, this method is somewhat limited as it already assumes knowledge of each term involved in the dynamics except for how they are scaled. The results of PINN are recorded in Table \ref{Tab:Maziar}. 
We found that even by increasing the number of parameters and training time, we are not guaranteed a good result.
\end{example}

\begin{table}[t]
\label{Tab:Maziar}
\vskip 0.15in
\begin{center}
\begin{small}
\begin{sc}
\begin{tabular}{lcccr}
\toprule
$\veps$ & Method & learned $g$-equation & Error Type-I/ II\\
\midrule
1/ 16  & PINN & $\partial_{t}g=-( 16 ^2-{\color{red} 1.253 \cdot 10^{ 1 }}) g -( 16 -{\color{red} 1.252 \cdot 10 ^{ 1 }}) v\cdot \partial_{x}g$  & \\
& & $+ ( 16 -{\color{red} 1.366 \cdot 10^{ 2 }})\langle v\partial_{x}g\rangle -( 16^2 -{\color{red} 3.097 \cdot 10^{ 1 }})v\cdot\partial_{x}\rho+\cdots$ & { 64.59 \%}\,/\,237.24\%\\

1/ 32  & PINN & $\partial_{t}g=-( 32 ^2-{\color{red} 4.990 \cdot 10^{ 1 }}) g -( 32 -{\color{red} 0.538 \cdot 10 ^{ 2 }}) v\cdot \partial_{x}g$  & \\
& & $+ ( 32 -{\color{red} 1.243 \cdot 10^{ 2 }})\langle v\partial_{x}g\rangle -( 32^2 -{\color{red} 1.256 \cdot 10^{ 2 }})v\cdot\partial_{x}\rho+\cdots$ & { 87.02 \%}\,/\,143.42\%\\

1/ 64  & PINN & $\partial_{t}g=-( 64 ^2-{\color{red} 0.700 \cdot 10^{ 2 }}) g -( 64 -{\color{red} 0.702 \cdot 10 ^{ 2 }}) v\cdot \partial_{x}g$  & \\
& & $+ ( 64 -{\color{red} 1.221 \cdot 10^{ 2 }})\langle v\partial_{x}g\rangle -( 4096 -{\color{red} 1.180 \cdot 10^{ 2 }})v\cdot\partial_{x}\rho+\cdots$ & { 95.72 \%}\,/\,76.26\%\\

1/ 128  & PINN & $\partial_{t}g=-( 128 ^2-{\color{red} 0.542 \cdot 10^{ 2 }}) g -( 128 -{\color{red} 0.543 \cdot 10 ^{ 2 }}) v\cdot \partial_{x}g$  & \\
& & $+ ( 128 -{\color{red} 1.223 \cdot 10^{ 2 }})\langle v\partial_{x}g\rangle -( 16384 -{\color{red} 1.221 \cdot 10^{ 2 }})v\cdot\partial_{x}\rho+\cdots$ & { 98.93 \%}\,/\,34.76\%\\

1/ 256  & PINN & $\partial_{t}g=-( 256 ^2-{\color{red}6.544 \cdot 10^{ 4 }}) g -( 256 -{\color{red}2.534 \cdot 10^{ 2 }}) v\cdot \partial_{x}g$  & \\

& & $+ ( 256 -{\color{red} 2.534 \cdot 10^{ 2 }})\langle v\partial_{x}g\rangle -( 65536 -{\color{red} 6.544 \cdot 10^{ 4 }})v\cdot\partial_{x}\rho+\cdots$ & { 99.86 \%}\,/\,99.42\%\\

\bottomrule
\end{tabular}
\end{sc}
\end{small}
 \end{center}\caption{Learned $g$-equation using PINN in \cite{DBLP:journals/corr/abs-1711-10561,DBLP:journals/corr/abs-1711-10566}.}
\vskip -0.1in
\end{table}

\begin{example}\label{compare4}
\textbf{Comparison with Conventional Methods. (Part 4)}
We now compare our results with the multiscale hierarchical deep learning (MS-HDL) approach proposed in \cite{Liu2020HierarchicalDL}. The approach in \cite{Liu2020HierarchicalDL} is to train separate feed-forward neural networks $\bd{F}_{j}(\bd{x},\Delta t_{j})$ for different time scales $\Delta t_{j}$:
\begin{equation}\label{eq:mshdlfwrd}
\bd{x}_{t+\Delta t_{j}}=\bd{x}_{t}+\bd{F}_{j}(\bd{x},\Delta t_{j}).
\end{equation}
For example, $\Delta t_{j}$ could be set to slow, medium, and fast scales by setting $\Delta t_{j}=\dfrac{\Delta t}{\veps^{j}}$ for some fixed $\veps$ and $j=0,1,2$. Unfortunately, \cite{Liu2020HierarchicalDL} does not provide a method for determining operators involved for each $\bd{F}_{j}(\bd{x},\Delta t_{j})$. Since we are interested in discovering the dynamics, we fit the $\bd{F}_{j}$ using the same 18 terms (denoted by $\mathcal{A}_{i}(v,x,t)$ for $i=1,2,cdots,18$) as in Example \ref{compare1}:
\begin{equation}
\bd{F}_{j}(\bd{x},t_{n}):=\sum_{i=1}^{18}\lambda_{i,j}\mathcal{A}_{i}(v,x,t_{n}).
\end{equation}
As suggested in Equation \eqref{eq:mshdlfwrd}, we propagate data using the forward Euler scheme. Thus, the $\lambda_{i,j}$ are determined using the loss in Equation \eqref{eq:fwrkappa}.

To be clear, Equation \eqref{eq:mshdlfwrd} is used to determine the dynamics of each map $\bd{F}_{j}(\bd{x},t)$ separately. Thus, the desired equations the MS-HDL would like to uncover are:
\begin{equation}\label{eq:Gsplit}
\begin{split}
\partial_{t}g_{fast}&=-\dfrac{\sigma^{A}}{\veps^{2}}g_{fast}-\dfrac{1}{\veps^{2}}v\partial_{x}\rho,\\
\partial_{t}g_{medium}&=-\dfrac{1}{\veps}(v\partial_{x}g_{medium}-\langle v\partial_{x}g_{medium}\rangle),\\
\partial_{t}g_{slow}&=-\sigma^{A}g_{slow}.
\end{split}
\end{equation}

We use $\sigma^{S}=1$, $\sigma^{A}=0$, and $G(x)=0$ to produce the data so that only fast and medium scales are present. For the MS-HDL, we choose $\Delta t_{j}=\frac{\Delta t}{\veps^{j}}$, $j=0,1,2$ with the correct value of $\veps$. In \cite{Liu2020HierarchicalDL}, the authors suggest gathering data for each time scale:
\begin{equation}\label{eq:msgather}
\begin{split}
g_{fast}(v,x,t_{n})&=g(v,x,n\Delta t_{2}), \hspace{1cm}n=0,1,2,...,N_{fast},\\
g_{medium}(v,x,t_{n})&=g(v,x,n\Delta t_{1}), \hspace{1cm}n=0,1,2,...,N_{medium},\\
g_{slow}(v,x,t_{n})&=g(v,x,n\Delta t_{0}), \hspace{1cm}n=0,1,2,...,N_{slow},
\end{split}
\end{equation}
i.e., the coarseness of the time grid determines the time scales. Of course, gathering data as in \eqref{eq:msgather} can be a problem. Namely, \eqref{eq:msgather} is only an approximation to the dynamics of \eqref{eq:Gsplit}. Thus, for our numerical example, we made the extra effort to perfectly split the data into different orders. In practice, it may be difficult to accurately split the data into different orders. For our DC-RNN algorithm, we do not need to split the data. The data for the DC-RNN is collected by:
\begin{equation}
g(v,x,t_{n})=g(v,x,n\Delta t)\hspace{1cm}n=0,1,2,...,N_{t}.
\end{equation}
Thus, one reason to prefer using DC-RNN over the MS-HDL is that one does not need to make the extra effort to split the data into different orders. Also, in the DC-RNN method we do not have to choose $\Delta t_{j}$ before hand, the DC-RNN algorithm learns appropriate time scales via Equation \eqref{eq:nonparallel} in an automatic manner. We compare our DC-RNN method with the MS-HDL method in Table \ref{Tab:KutzMethod}.
\end{example}

\begin{table}[t]
\label{Tab:KutzMethod}
\vskip 0.15in
\begin{center}
\begin{small}
\begin{sc}
\begin{tabular}{lccr}
\toprule
$\veps$ & Method & Learned $g$-equation & Error Type-I/ II\\
\midrule
1/ 16  & MS-HDL & $\partial_{t}g=-( 16 ^2+{\color{red} 3.663}) g -( 16 +{\color{red} 4.737 \cdot 10 ^{ -1 }}) v\cdot 
\partial_{x}g$  & \\
& & $+ ( 16 +{\color{red} 0.741})\langle v\partial_{x}g\rangle -( 256 -{\color{red} 3.315 \cdot 10^{ -1 }})v\cdot\partial_{x}\rho+\cdots$ & { 0.96 \%}\,/\,2.29\%\\

1/ 32  & MS-HDL & $\partial_{t}g=-( 32 ^2+{\color{red} 2.521}) g -( 32 +{\color{red} 4.927 \cdot 10 ^{ -1 }}) v\cdot \partial_{x}g$  & \\
& & $+ ( 32 +{\color{red} 1.209})\langle v\partial_{x}g\rangle -( 1024 -{\color{red} 0.642 \cdot 10^{ 1 }})v\cdot\partial_{x}\rho+\cdots$ & { 0.50 \%}\,/1.55\%\\

1/ 64  & MS-HDL & $\partial_{t}g=-( 64 ^2-{\color{red} 0.954 \cdot 10^{ 2 }}) g -( 64 -{\color{red} 1.617}) v\cdot \partial_{x}g$  & \\
& & $+ ( 64 +{\color{red} 2.100 \cdot 10^{ -1 }})\langle v\partial_{x}g\rangle -( 4096 -{\color{red} 1.110 \cdot 10^{ 2 }})v\cdot\partial_{x}\rho+\cdots$ & { 2.50 \%}\,/\,1.97\%\\

1/ 128  & MS-HDL & $\partial_{t}g=-( 128 ^2-{\color{red} 0.995 \cdot 10^{ 3 }}) g -( 128 -{\color{red} 1.231 \cdot 10 ^{ 1 }}) v\cdot \partial_{x}g$  & \\
& & $+ ( 128 -{\color{red} 0.786 \cdot 10^{ 1 }})\langle v\partial_{x}g\rangle -( 16384 -{\color{red} 1.050 \cdot 10^{ 3 }})v\cdot\partial_{x}\rho+\cdots$ & { 6.26 \%}\,/\,7.06\%\\

1/ 256  & MS-HDL & $\partial_{t}g=-( 256 ^2-{\color{red} 1.470 \cdot 10^{ 4 }}) g -( 256 -{\color{red} 0.829 \cdot 10 ^{ 2 }}) v\cdot \partial_{x}g$  & \\
& & $+ ( 256 -{\color{red} 0.572 \cdot 10^{ 2 }})\langle v\partial_{x}g\rangle -( 65536 -{\color{red} 1.470 \cdot 10^{ 4 }})v\cdot\partial_{x}\rho+\cdots$ & { 22.45 \%}\,/\,24.90\%\\

\bottomrule
\end{tabular}
\end{sc}
\end{small}
\end{center}\caption{Learned $g$-equation using the multiscale hierarchical deep learning method (MS-HDL) in \cite{Liu2020HierarchicalDL}.}
\vskip -0.1in
\end{table}

\section{Conclusion}
\label{sec:conclusion}
We propose a deep learning algorithm capable of learning time-dependent multiscale and nonlocal partial differential equations (PDEs) from data. The key to achieving our goal is to construct a Densely Connected Recurring Neural Network (DC-RNN) that accounts for potential multiscale and nonlocal structures in the data. The DC-RNN is a symbolic network with relationship among the symbols given by high-order IMEX schemes used to target dynamics of multiscale kinetic equations. Incorporated into the training of the network are physics-aware constraints and multiscale ansatz. Through various numerical experiements, we verify that our DC-RNN accurately and efficiently recovers multiscale PDEs which the data satisfies. As a byproduct, our DC-RNN determines appropriate multiscale parameters and can potentially discover lower dimensional representations for kinetic equations.

\section*{Acknowledgments} H.~Yang was partially supported by the US National Science Foundation under award DMS-1945029. 

\section{Appendix}

Here we present details on how to define a loss function which makes use of high-order IMEX schemes to fit data to Equations \eqref{eq:rho} and \eqref{eq:GGG}.

\subsection{High-order IMEX Runge-Kutta fitting}
\label{imex1}


A K-stage IMEX Runge-Kutta scheme is given by:
\begin{equation}\label{eq:intermediateg}
\begin{split}
g^{(i)}&=g^{n}-\Delta t\sum_{j=1}^{i-1}\wt{a}_{i,j}\left(\dfrac{1}{\veps}(I-\langle\rangle)(v\partial_{x}g^{(j)})+\dfrac{1}{\veps^{2}}v\partial_{x}\rho^{(j)}+\sigma^{A}g^{(j)}\right)-\Delta t\sum_{j=1}^{i}a_{i,j}\left(\dfrac{\sigma^{S}}{\veps^{2}}g^{(j)}\right),
\end{split}
\end{equation}
\begin{equation}\label{eq:intermediaterho}
\rho^{(i)}=\rho^{n}-\Delta t \sum_{j=1}^{i-1}\wt{a}_{i,j}(\sigma^{A}\rho^{(j)}-G)-\Delta t\sum_{j=1}^{i}a_{i,j}\partial_{x}\langle v g^{(j)}\rangle,
\end{equation}
\begin{equation}\label{eq:gnplus1}
\begin{split}
g^{n+1}&=g^{n}-\Delta t \sum_{i=1}^{K}\wt{w}_{i}\left(\dfrac{1}{\veps}(I-\langle \rangle)(v\partial_{x}g^{(i)})+\dfrac{1}{\veps^{2}}v\partial_{x}\rho^{(i)}+\sigma^{A}g^{(i)}\right)-\Delta t\sum_{j=1}^{K}w_{i}\left(\dfrac{\sigma^{S}}{\veps^{2}}g^{(j)}\right),
\end{split}
\end{equation}
\begin{equation}\label{eq:rhonplus1}
\rho^{n+1}=\rho^{n}-\Delta t \sum_{i=1}^{K}\wt{w}_{i}(\sigma^{A}\rho^{(i)}-G)-\Delta t\sum_{i=1}^{K}w_{i}\partial_{x}\langle v g^{(i)}\rangle.
\end{equation}

Equations \eqref{eq:intermediateg} and \eqref{eq:intermediaterho} are intermediate stages and Equations \eqref{eq:gnplus1} and \eqref{eq:rhonplus1} are the approximate solution at the next time step. Here $\wt{A}=(\wt{a}_{i,j})$ with $\wt{a}_{i,j}=0$ for $j\geq i$ and $A=(a_{i,j})$ with $a_{i,j}=0$ for $j>i$ are $K\times K$ matrices. Along with the coefficient vectors $\wt{\bd{w}}=(\wt{w}_{1},\cdots,\wt{w}_{K})^{T}$, $\bd{w}=(w_{1},\cdots,w_{K})^{T}$,
 they can be represented by a double Butcher tableau:
\begin{center}
$\renewcommand\arraystretch{1.2}
\begin{array}
{c|cccc}
\wt{\bd{c}}& \wt{A}\\
\hline
& \wt{\bd{w}}^{T}
\end{array}\hspace{1cm} \text{and} \hspace{1cm}
\renewcommand\arraystretch{1.2}\begin{array}{c|cccc}
\bd{c}& A\\
\hline
& \bd{w}^{T}
\end{array}$,
\end{center}
where the vectors $\wt{\bd{c}}=(\wt{c}_{1},\cdots,\wt{c}_{K})^{T}$ and $\bd{c}=(c_{1},\cdots,c_{K})^{T})$ are defined as:
\begin{equation}
\wt{c}_{i}=\sum_{j=1}^{i-1}\wt{a}_{i,j}\hspace{1cm}\text{and}\hspace{1cm} c_{i}=\sum_{j=1}^{i-1}a_{i,j}.
\end{equation}

For convenience, we provide the tableau for the ARS(2,2,2) scheme:
\[
\renewcommand\arraystretch{1.2}
\begin{array}
{c|ccc}
0 & 0 & 0 & 0\\
\gamma & \gamma & 0 & 0 \\
1& \delta & 1-\delta & 0\\
\hline
& \delta & 1-\delta & 0
\end{array}\hspace{1cm} \text{and} \hspace{1cm}
\renewcommand\arraystretch{1.2}
\begin{array}
{c|ccc}
0 & 0 & 0 & 0\\
\gamma & 0 & \gamma & 0 \\
1& 0 & 1-\gamma & \gamma\\
\hline
& 0 & 1-\gamma & \gamma
\end{array},
\]
where $\gamma=1-\dfrac{\sqrt{2}}{2}$ and $\delta=1-\dfrac{1}{2\gamma}$. 

\comment{This is a second-order IMEX scheme is tested in the Numerical.}

The loss function based on this fitting scheme is defined by:
\begin{equation}
L=\dfrac{1}{Nt-1}\sum_{n=1}^{Nt-1}||\mathcal{K}_{g}^{n}||+||\mathcal{K}_{\rho}^{n}||
\end{equation}
with
\begin{equation}
\mathcal{K}_{g}^{n}=\mathcal{K}_{g}^{n}(\{g(v,x,t_{n}),g(v,x,t_{n+1})\}), \quad \mathcal{K}_{\rho}^{n}=\mathcal{K}_{\rho}^{n}(\{\rho(x,t_{n}),\rho(x,t_{n+1})\}),
\end{equation}
where
\begin{equation}
\begin{split}
\mathcal{K}_{g}^{n}&:=g(v,x,t_{n+1})-g(v,x,t_{n})+\Delta t\left(\sum_{i=1}^{K}\sigma^{A}(x)\wt{w}_{i}g^{(i)}+\dfrac{\sigma^{S}(x)}{\veps^{2}}w_{i}g^{(i)}\right)+\Delta t \sum_{i=1}^{K}\wt{w}_{i}\left(\mathcal{F}_{1}(g^{(i)}(v,x),\rho^{(i)}(x))\right)
\end{split},
\end{equation}

\begin{equation}
\begin{split}
\mathcal{K}_{\rho}^{n}&:=\rho(x,t_{n+1})-\rho(x,t_{n})+\Delta t\sum_{i=1}^{K}\wt{w}_{i}(\sigma^{A}(x)\rho^{(i)}-G(x))+\Delta t\sum_{i=1}^{K}w_{i}\left(\mathcal{F}_{2}(g^{(i)}(v,x),\rho^{(i)}(x))\right).
\end{split}
\end{equation}

The operators $\mathcal{F}_{1}(g,\rho)$, $\mathcal{F}_{2}(g,\rho)$ are given by \eqref{eq:vectorizedansatz} and are generated by the RNN in Equation \eqref{eq:eq210}. The intermediate stages are given by:

\begin{equation}
\begin{split}
g^{(i)}&=g(v,x,t_{n})-\Delta t\sum_{j=1}^{i}a_{i,j}\dfrac{\sigma^{S}(x)}{\veps^{2}}g^{(j)}-\Delta t\sum_{j=1}^{i-1}\wt{a}_{i,j}\left(\mathcal{F}_{1}(g^{(j)},\rho^{(j)})\right),
\end{split}
\end{equation}

\begin{equation}
\begin{split}
\rho^{(i)}&=\rho(x,t_{n})-\Delta t\sum_{j=1}^{i-1}\wt{a}_{i,j}(\sigma^{A}(x)\rho^{(j)}-G)-\Delta t\sum_{j=1}^{i}\wt{a}_{i,j}\left(\mathcal{F}_{2}(g^{(j)},\rho^{(j)})\right).
\end{split}
\end{equation}

We note that $\sigma^{A}(x)$, $\sigma^{S}(x)$, and $G(x)$ do not need to be assumed known. These functions can be part of the fitting process by replacing them with feed-forward neural nets, say.

\subsection{High-order IMEX-BDF fitting}
\label{imex2}

Another way to go higher order in time is through the IMEX-BDF scheme:
\begin{equation}
\begin{split}
\sum_{i=0}^{q}\alpha_{i}g^{n+i}&+\Delta t\sum_{i=0}^{q-1}\gamma_{i}\left(\dfrac{1}{\veps}(I-\langle\rangle)(v\partial_{x}g^{n+i})\right.\\
&\left.+\dfrac{1}{\veps^{2}}v\partial_{x}\rho^{n+i}+\sigma^{A}g^{n+i}\right)+\beta\Delta t \left(\dfrac{\sigma^{S}}{\veps^{2}}g^{n+q}\right)=0,
\end{split}
\end{equation}
and
\begin{equation}
\sum_{i=0}^{q}\alpha_{i}\rho^{n+i}+\Delta t\sum_{i=0}^{q-1}\gamma_{i}(\sigma^{A}\rho^{n+i}-G)+\beta\Delta t \partial_{x}\langle v g^{n+q}\rangle=0.
\end{equation}

We display some coefficients $\alpha=(\alpha_{0},\cdots,\alpha_{q})$, $\gamma=(\gamma_{0},\cdots,\gamma_{q-1})$, and $\beta$ for the above scheme in Table \ref{Tab:table-1}.

\begin{table}[t]
\caption{}
\label{Tab:table-1}
\begin{center}
\begin{small}
\begin{sc}
\begin{tabular}{lccr}
\toprule
$q$ & $\alpha$ & $\gamma$ & $\beta$ \\
\midrule
1 & $(-1,1)$ & $1$ & $1$\\
\vspace{.1cm}
2 & $(\frac{1}{3},-\frac{4}{3},1)$ & $(-\frac{2}{3},\frac{4}{3})$ & $\frac{2}{3}$\\
\vspace{.1cm}
3 & $(-\frac{2}{11},\frac{9}{11},-\frac{18}{11},1)$ & $(\frac{6}{11},-\frac{18}{11},\frac{18}{11})$ & $\frac{6}{11}$\\
\vspace{.1cm}
4 & $(\frac{3}{25},-\frac{16}{25},\frac{36}{25},-\frac{48}{25},1)$ & $(-\frac{12}{25},\frac{48}{25},-\frac{72}{25},\frac{48}{25})$ & $\frac{12}{25}$\\
\bottomrule
\end{tabular}
\end{sc}
\end{small}
\end{center}
\end{table}

The loss function for the fitting scheme based on the IMEX-BDF scheme is defined by:
\begin{equation}
L=\dfrac{1}{N_{t}-q}\sum_{n=1}^{Nt-q}||\mathcal{K}^{n}(D;\bd{\theta})||
\end{equation}
with
\begin{equation}
D=\{u(x,t_{n}),u(x,t_{n+1})\cdots,u(x,t_{n+q})\}.
\end{equation}

For the $g$ equation $\mathcal{K}_{g}^{n}$ is given by:
\begin{equation}
\begin{split}
\mathcal{K}_{g}^{n}=&\sum_{i=0}^{q}\alpha_{i}g^{n+i}-\beta\Delta t \dfrac{\sigma^{S}(x)}{\veps^{2}}g^{n+q}-\Delta t \sum_{i=0}^{q-1}\sigma^{A}(x)g^{n+i}+\Delta t\sum_{i=0}^{q-1}\gamma_{i}\left(\mathcal{F}_{1}(g(v,x,t_{n+i}),\rho(x,t_{n+i}))\right).
\end{split}
\end{equation}
The operator $\mathcal{F}_{1}(g,\rho)$ is given by \eqref{eq:vectorizedansatz} and is generated by the RNN in Equation \eqref{eq:eq210}.

For the $\rho$ equation $\mathcal{K}_{g}^{n}$ is given by:
\begin{equation}
\begin{split}
\mathcal{K}_{\rho}^{n}=&\sum_{i=0}^{q}\alpha_{i}g^{n+i}+\Delta t \sum_{i=0}^{q-1}\gamma_{i}\left(\sigma^{A}\rho^{n+i}-G\right)-\beta\Delta t \left(\mathcal{F}_{2}(g(v,x,t_{n+q}),\rho(x,t_{n+q}))\right).
\end{split}
\end{equation}

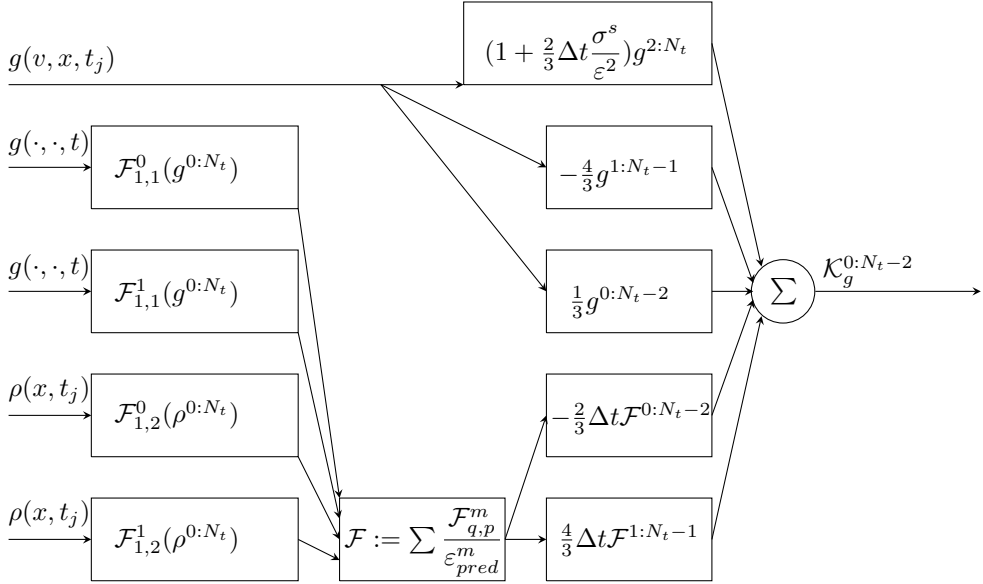
\begin{figure}\label{NNfig2}
\begin{tikzpicture}[x=1.1cm, y=1.1cm, >=stealth]

\draw [black, ->          ] (0,4.5) -- (5.5,4.5);
\put(0,148){$g(v,x,t_{j})$}

\draw [black, ->          ] (4.5,4.5) -- (6.5,3.5);

\draw [black, ->          ] (4.5,4.5) -- (6.5,2.0);

\draw [black, ->          ] (6.0,-1.0) -- (6.5,0.5);

\draw [black, ->          ] (6.0,-1.0) -- (6.5,-1.0);

\draw [black, ->          ] (0,3.5) -- (1,3.5);
\put(0,116){$g(\cdot,\cdot,t)$}

\put(40,106){$\mathcal{F}_{1,1}^{0}(g^{0:N_{t}})$}
\draw[black] (1,3) -- (1,4) -- (3.5,4) -- (3.5,3) -- (1,3);

\draw [black, ->          ] (0,2) -- (1,2);
\put(0,69){$g(\cdot,\cdot,t)$}

\put(40,59){$\mathcal{F}_{1,1}^{1}(g^{0:N_{t}})$}
\draw[black] (1,1.5) -- (1,2.5) -- (3.5,2.5) -- (3.5,1.5) -- (1,1.5);

\draw [black, ->          ] (0,0.5) -- (1,0.5);
\put(0,23){$\rho(x,t_{j})$}

\put(40,13){$\mathcal{F}_{1,2}^{0}(\rho^{0:N_{t}})$}
\draw[black] (1,0) -- (1,1) -- (3.5,1) -- (3.5,0) -- (1,0);

\draw [black, ->          ] (0,-1) -- (1,-1);
\put(0,-24){$\rho(x,t_{j})$}

\put(40,-34){$\mathcal{F}_{1,2}^{1}(\rho^{0:N_{t}})$}
\draw[black] (1,-1.5) -- (1,-0.5) -- (3.5,-0.5) -- (3.5,-1.5) -- (1,-1.5);

\put(127,-34){$\mathcal{F}:=\sum\dfrac{\mathcal{F}_{q,p}^{m}}{\veps_{pred}^{m}}$}
\draw (4,-1.5) -- (4,-0.5) -- (6,-0.5) -- (6,-1.5) -- (4,-1.5);

\draw [black, ->  ] (3.5,3) -- (4,-0.5);

\draw [black, ->  ] (3.5,1.5) -- (4,-0.75);

\draw [black, ->  ] (3.5,0) -- (4,-1);

\draw [black, ->  ] (3.5,-1.0) -- (4,-1.25);

\put(180,150){$(1+\frac{2}{3}\Delta t \dfrac{\sigma^{s}}{\veps^{2}})g^{2:N_t}$}
\draw (5.5,4.5) -- (5.5,5.5) -- (8.5,5.5) -- (8.5,4.5) -- (5.5,4.5);

\put(208,104){$-\frac{4}{3} g^{1:N_t-1}$}
\draw (6.5,3) -- (6.5,4) -- (8.5,4) -- (8.5,3) -- (6.5,3);

\put(212,56){$\frac{1}{3} g^{0:N_t-2}$}
\draw (6.5,1.5) -- (6.5,2.5) -- (8.5,2.5) -- (8.5,1.5) -- (6.5,1.5);

\put(208,-34){$\frac{4}{3}\Delta t \mathcal{F}^{1:N_t-1}$}
\draw (6.5,-1.5) -- (6.5,-0.5) -- (8.5,-0.5) -- (8.5,-1.5) -- (6.5,-1.5);

\put(205,12){$-\frac{2}{3}\Delta t \mathcal{F}^{0:N_t-2}$}
\draw (6.5,0) -- (6.5,1) -- (8.5,1) -- (8.5,0) -- (6.5,0);

\draw [black, ->  ] (8.5,-1.0) -- (9.1,1.7);

\draw [black, ->  ] (8.5,0.5) -- (9,1.9);

\draw [black, ->  ] (8.5,2) -- (9,2);

\draw [black, ->  ] (8.5,3.5) -- (9,2.1);

\draw [black, ->  ] (8.5,5) -- (9.1,2.3);

\draw[] (9.36,2) circle (12pt);
\put(288,60){$\sum$}

\put(308,68){$\mathcal{K}_{g}^{0:N_t-2}$}
\draw [black, ->  ] (9.75,2) -- (11.75,2);

\end{tikzpicture}
\caption{Example DC-RNN based on IMEX-BDF-2 scheme for predicting the $g$-equation. The inputs are $\rho(t)$, and $g(t)$. The dictionary contains order $\mathcal{O}(1)$ and $\mathcal{O}(\veps)$ operators. These operators are generated by the RNNs of orders $\veps^{-m}$ $m=0,1$. The output $\mathcal{K}_{g}^{0:N_t-2}$ is to be minimized with respect to a chosen norm.}
\end{figure}

Again, $\sigma^{A}(x)$, $\sigma^{S}(x)$, and $G(x)$ can be learned by including them in the fitting process. We display in Figure \ref{NNfig2} a DC-RNN for determining the equation satisfied by $g(v,x,t)$ based on the IMEX-BDF2 scheme.

\bibliographystyle{siamplain}
\bibliography{ref}
\end{document}